\documentclass[preprint,dvipdfmx]{elsarticle}
\usepackage{bm,amsmath,amssymb,amsthm,revsymb,ascmac}
\newtheoremstyle{mydefinition}{}{}{}{}{\bfseries}{}{\newline}{\thmname{#1}\thmnumber{\ #2}\thmnote{\quad#3}.}
\theoremstyle{mydefinition}
\newtheorem{Thm}{Theorem} 
\newtheorem{Cor}[Thm]{Corollary} 
\newtheorem{Prop}[Thm]{Proposition}
\newtheorem{Lem}{Lemma}
\newtheorem*{Lem*}{Lemma}

\makeatletter
\newcommand{\opnorm}{\@ifstar\@opnorms\@opnorm}
\newcommand{\@opnorms}[1]{%
\left|\mkern-1.5mu\left|\mkern-1.5mu\left|
#1
\right|\mkern-1.5mu\right|\mkern-1.5mu\right|
}
\newcommand{\@opnorm}[2][]{%
\mathopen{#1|\mkern-1.5mu#1|\mkern-1.5mu#1|}
#2
\mathclose{#1|\mkern-1.5mu#1|\mkern-1.5mu#1|}
}
\newenvironment{Prf}[1][]{\par%%
\pushQED{\qed}%%
\normalfont \topsep6\p@\@plus6\p@\relax%%
\trivlist%%
\item[\hskip\labelsep%%
\@ifempty{#1}%%
{{\itshape Proof}}%%
{{\itshape Proof of\,\,\,#1}}
\@addpunct{.}%%
]%\ignorespaces%%
\mbox{}\linebreak%%
}{%%
\popQED\endtrivlist\@endpefalse%%
}%%
\makeatother
\makeatletter

\@addtoreset{equation}{section}
\makeatother
\def\hsmash#1{\makebox[0pt][l]{#1}}
\begin{document}
%%%%%%%%%%%%%%%%%%%%%%%%%%%%%%%%%%%%%%%%%%%%%%%%%%
%%%%%%%%%%%%%%%%%%%%%%%%%%%%%%%%%%%%%%%%%%%%%%%%%%
%%%%%%%%%%   Front Matter   %%%%%%%%%%%%%%%%%%%%%%
%%%%%%%%%%%%%%%%%%%%%%%%%%%%%%%%%%%%%%%%%%%%%%%%%%
%%%%%%%%%%%%%%%%%%%%%%%%%%%%%%%%%%%%%%%%%%%%%%%%%%
\begin{frontmatter}
\title{Local well-posedness of the complex Ginzburg-Landau equation in bounded domains}

\author[add1]{Takanori Kuroda}
\ead{1d\_est\_quod\_est@ruri.waseda.jp}
\author[add2,fn2]{Mitsuharu \^Otani}
\ead{otani@waseda.jp}
\fntext[fn2]{Partly supported by the Grant-in-Aid for Scientific Research, \# 15K13451, the Ministry of Education, Culture, Sports, Science and Technology, Japan.}

\address[add1]{Major in Pure and Applied Physics,\\ Graduate School of Advanced Science and Engineering, \\ Waseda University, 3-4-1 Okubo Shinjuku-ku, Tokyo, 169-8555, JAPAN}
\address[add2]{Department of Applied Physics, School of Science and Engineering, \\ Waseda University, 3-4-1 Okubo Shinjuku-ku, Tokyo, 169-8555, JAPAN}
%%%%%%%%%%%%%%%%%%%%%%%%%%%%%%%%%%%%%%%%%%%%%%%%%%
%%%%%%%%%%   Abstract   %%%%%%%%%%%%%%%%%%%%%%%%%%
%%%%%%%%%%%%%%%%%%%%%%%%%%%%%%%%%%%%%%%%%%%%%%%%%%
\begin{abstract}
In this paper, we are concerned with the local well-posedness of the initial-boundary value problem for complex Ginzburg-Landau (CGL) equations in bounded domains.
There are many studies for the case where the real part of its nonlinear term plays as dissipation.
This dissipative case is intensively studied and it is shown that (CGL) admits a global solution when parameters appearing in (CGL) belong to the so-called CGL-region.
This paper deals with the non-dissipative case.
We regard (CGL) as a parabolic equation perturbed by monotone and non-monotone perturbations and follows the basic strategy developed in \^Otani (1982) to show  the local well-posedness of (CGL) and the existence of small global solutions provided that the nonlinearity is Sobolev subcritical.
\end{abstract}
%%%%%%%%%%%%%%%%%%%%%%%%%%%%%%%%%%%%%%%%%%%%%%%%%%
%%%%%%%%%%   Key Word   %%%%%%%%%%%%%%%%%%%%%%%%%%
%%%%%%%%%%%%%%%%%%%%%%%%%%%%%%%%%%%%%%%%%%%%%%%%%%
\end{frontmatter}
%%%%%%%%%%%%%%%%%%%%%%%%%%%%%%%%%%%%%%%%%%%%%%%%%%
%%%%%%%%%%%%%%%%%%%%%%%%%%%%%%%%%%%%%%%%%%%%%%%%%%
%%%%%%%%%%   Introduction   %%%%%%%%%%%%%%%%%%%%%%
%%%%%%%%%%%%%%%%%%%%%%%%%%%%%%%%%%%%%%%%%%%%%%%%%%
%%%%%%%%%%%%%%%%%%%%%%%%%%%%%%%%%%%%%%%%%%%%%%%%%%
\section{Introduction}
In this paper, we consider the following Cauchy problem for the complex Ginzburg-Landau equation.
\begin{equation}
\tag*{(CGL)}
\left\{
\begin{aligned}
&\frac{\partial u}{\partial t}(t, x) \!-\! (\lambda \!+\! i\alpha)\Delta u \!-\! (\kappa \!+\! i\beta)|u|^{q-2}u \!-\! \gamma u \!=\! f(t, x),\\
&\hspace{50mm}\mbox{in}\ (t, x) \in [0, T] \times \Omega,\\
&\hsmash{\(u(t, x) = 0\)}\hspace{50mm}\mbox{on}\ (t, x) \in [0, T] \times \partial\Omega,\\
&\hsmash{\(u(0, x) = u_0(x)\)}\hspace{50mm}\mbox{in}\ x \in \Omega,
\end{aligned}
\right.
\end{equation}
where \(\Omega\) is a bounded domain in \(\mathbb{R}^N\) with smooth boundary \(\partial\Omega\); 
\(\lambda, \kappa\) are positive parameters, \(\alpha, \beta, \gamma\) are real parameters and \(i = \sqrt{-1}\) denotes the imaginary unit; 
\(q \geq 2\) is a given number; 
\(f: [0, T] \to \mathbb{C}\) is a given external force defined on an interval \([0, T]\) with \(T > 0\).
Our unknown function \(u: [0, T] \to \mathbb{C}\), which takes values in the complex numbers, represents an order parameter.

%%%%%%%%%%%%%%%%%%%%%%%%%%%%%%%%%%%%%%%%%%%%%%%%%%
%%%%%%%%%%   Physical background   %%%%%%%%%%%%%%%
%%%%%%%%%%%%%%%%%%%%%%%%%%%%%%%%%%%%%%%%%%%%%%%%%%
This equation is originally introduced by Ginzburg \& Landau \cite{GL1} in order to give a mathematical model of superconductivities.
In their theory, the Ginzburg-Landau theory for superconductivities, the physical quantity \(u\) is only a parameter which describes the randomness of a certain physical state.
In the sequel, the theory has been applied to more general phenomena related to phase change (Nishiura \cite{N1}) or pattern formulation (Cross \& Hohenberg \cite{CH1}).

%%%%%%%%%%%%%%%%%%%%%%%%%%%%%%%%%%%%%%%%%%%%%%%%%%
%%%%%%%%%%   Math approach k < 0   %%%%%%%%%%%%%%%
%%%%%%%%%%%%%%%%%%%%%%%%%%%%%%%%%%%%%%%%%%%%%%%%%%
There have been many contribution to the case where \(\kappa < 0\): 
Temam \cite{T1} showed the existence of the unique weak solution using Galerkin method for the case where \(N = 1, 2\) and \(q = 4\); 
Yang \cite{Y1} got mild solutions in terms of the semi-group \(\left\{e^{(\lambda+i\alpha)t \Delta}; t \geq 0 \right\}\) for the case where \(N = 1, 2, 3\) and \(q \leq 2 + \frac{4}{N}\); 
Okazawa \& Yokota \cite{OY1} developed the maximal monotone operator theory in Hilbert spaces over the complex field and proved that strong solutions exist in bounded \(\Omega\) for the case where \(\left(\frac{\alpha}{\lambda}, \frac{\beta}{\kappa}\right)\) belongs to the so-called ``CGL-region'' and \(u_0 \in {\rm H}^1_0 \cap {\rm L}^q\) as well as the case where \(|\beta|/\kappa \leq 2\sqrt{q-1}/{(q-2)}\) and \(u_0 \in {\rm L}^2\); 
in the sequel, Kuroda, \^Otani \& Shimizu \cite{KOS1} improved this result by showing the existence of global solutions in general domains with \(\left(\frac{\alpha}{\lambda}, \frac{\beta}{\kappa}\right)\) being in the ``CGL-region'' and \(u_0 \in {\rm L}^2\).
In \cite{KOS1}, the abstract theory of parabolic equations are used, where \(-\lambda\Delta u\) is regarded as its leading term, \(-i\alpha\Delta u\) as a monotone perturbation and \(-i\beta |u|^{q-2}u\) as a non-monotone perturbation.
We remark that when \(\kappa < 0\) the real part of our nonlinear term \(-\kappa|u|^{q-2}u\) becomes maximal monotone in \({\rm L}^2\) and \(-\lambda\Delta u -\kappa|u|^{q-2}u\) can be represented as a subdifferential operator.

%%%%%%%%%%%%%%%%%%%%%%%%%%%%%%%%%%%%%%%%%%%%%%%%%%
%%%%%%%%%%   Math approach k > 0   %%%%%%%%%%%%%%%
%%%%%%%%%%%%%%%%%%%%%%%%%%%%%%%%%%%%%%%%%%%%%%%%%%
On the other hand, there are few treatment for the case where \(\kappa > 0\) especially on the well-posedness: 
Cazenave et al. studied blow-up of solutions(\cite{CDW1, CDF1}) and the existence of standing wave solutions (\cite{CDW2}).
In their papers \cite{CDW1, CDF1} they made restriction on parameters \(\lambda, \kappa, \alpha, \beta\) to be \(\frac{\alpha}{\lambda} = \frac{\beta}{\kappa}\).
They proved the existence of the unique local solution based on the semi-group theory for sufficiently smooth initial data in the whole space \(\mathbb{R}^N\) for any \(q > 2\).
The local well-posedness of (CGL) for general coefficients was treated in Shimotsuma, Yokota \& Yoshii \cite{SYY2}.
They showed the local existence of the unique solution in various kind of domains using the semi-group theory in \({\rm L}^p\) over the complex numbers under the assumption \(2 < q < 2 + \frac{2p}{N}\).
They also deduced the global extension of solutions by assuming \(\frac{|\alpha|}{\lambda} < \frac{2\sqrt{p - 1}}{p - 2}\) and a suitable condition on \(\gamma\).

%%%%%%%%%%%%%%%%%%%%%%%%%%%%%%%%%%%%%%%%%%%%%%%%%%
%%%%%%%%%%   Our approach   %%%%%%%%%%%%%%%%%%%%%%
%%%%%%%%%%%%%%%%%%%%%%%%%%%%%%%%%%%%%%%%%%%%%%%%%%
The main purpose of this paper is to discuss the local well-posedness of (CGL) in a bounded domain \(\Omega\) when \(u_0\) belongs to \({\rm H}_0^1(\Omega)\) and \(q\) is Sobolev subcritical, i.e.,
\[
2 \leq q < 2^* := 
\begin{cases}
+\infty & N = 1,2,\\
\frac{2N}{N-2} & N \geq 3.
\end{cases}
\]
To expect the local well-posedness of (CGL) under this situation is quite natural on the analogy of the theory of nonlinear parabolic equations.
In order to show that this conjecture holds true, we follow the basic strategy in \cite{KOS1}, i.e., we regard (CGL) as a parabolic equation with the principal part \(-\lambda \Delta u\) perturbed by the monotone perturbation \(-i \alpha\Delta u\) and the non-monotone perturbations \(-(\kappa + i\beta)|u|^{q-2}u\).
To cope with these perturbations, as for the monotone perturbation, we can use the standard argument from the maximal monotone operator theory.
As for the non-monotone perturbations, we rely on the Schauder-Tychonoff fixed point theory as in \cite{O2}.

%%%%%%%%%%%%%%%%%%%%%%%%%%%%%%%%%%%%%%%%%%%%%%%%%%
%%%%%%%%%%   Outline   %%%%%%%%%%%%%%%%%%%%%%%%%%%
%%%%%%%%%%%%%%%%%%%%%%%%%%%%%%%%%%%%%%%%%%%%%%%%%%
The outline of this paper is as follows.
In Section 2, we fix function spaces over the real numbers which are direct products of usual Lebesgue or Sobolev spaces over the real numbers.
Introducing suitable functionals and their subdifferentials on this space, we rewrite (CGL) in terms of  an evolution equation in a real Hilbert space governed by subdifferential operators with perturbations.
We conclude this section by stating our main results on the local well-posedness (Theorem \ref{lwpndbd}), an alternative on the maximal existence time of solutions (Theorem \ref{altbd}) and the existence of small global solutions (Theorem \ref{gebd}).

In Section 3, we discuss the solvability of some auxiliary equations (Proposition \ref{gwpAEh}), which is needed for the application of the Schauder-Tychonoff fixed point theorem.
In Section 4, we establish some a priori estimates for the solutions of auxiliary equations, by which we prove the existence part of Theorem \ref{lwpndbd}.
Theorem \ref{altbd} will be proved in Section 5 and the uniqueness part of Theorem \ref{lwpndbd} is shown in Section 6.
The last section is devoted to the proof of Theorem \ref{gebd}.

%%%%%%%%%%%%%%%%%%%%%%%%%%%%%%%%%%%%%%%%%%%%%%%%%%
%%%%%%%%%%%%%%%%%%%%%%%%%%%%%%%%%%%%%%%%%%%%%%%%%%
%%%%%%%%%%   Preliminaries   %%%%%%%%%%%%%%%%%%%%%
%%%%%%%%%%%%%%%%%%%%%%%%%%%%%%%%%%%%%%%%%%%%%%%%%%
%%%%%%%%%%%%%%%%%%%%%%%%%%%%%%%%%%%%%%%%%%%%%%%%%%
\section{Preliminaries}
In this section, we fix some function spaces to be used later and prepare notations in order to reduce (CGL) to an evolution equation in a certain product function space.
With these preparations, this section will be concluded by stating our main results of this paper.
We first introduce product function spaces made up of usual Lebesgue and Sobolev spaces over the real field using the following identification:
\[
\mathbb{C} \ni u_1 + iu_2 \mapsto U = (u_1, u_2) \in \mathbb{R}^2.
\]
These spaces are also Banach or Hilbert spaces with respect to the following norms or inner products. 
We set 
\begin{itemize}
\item[] $\bullet$ \ \(\mathbb{L}^r(\Omega) := {\rm L}^r(\Omega) \times {\rm L}^r(\Omega) \ni U = (u_1, u_2), V = (v_1, v_2)\),
\\[1mm]
 \quad with norm: \(|U|_{\mathbb{L}^r}^r = |u_1|_{{\rm L}^r}^r + |u_2|_{{\rm L}^r}^r\),
\\[1mm]
\quad and inner product (\(r = 2\)): \((U, V)_{\mathbb{L}^2} = (u_1, v_1)_{{\rm L}^2} + (u_2, v_2)_{{\rm L}^2}\).
\item[] $\bullet$ \ \(\mathbb{H}_0^1(\Omega) := {\rm H}_0^1(\Omega) \times {\rm H}_0^1(\Omega) \ni U = (u_1, u_2), V = (v_1, v_2)\),
\\[1mm]
\quad with inner product: \((U, V)_{\mathbb{H}^1} = (u_1, v_1)_{{\rm H}^1} + (u_2, v_2)_{{\rm H}^1}\).
\end{itemize}
In addition, \(\mathcal{H}^S\) denotes the space of functions with value in \(\mathbb{L}^2(\Omega)\) defined on \([0, S]\) with \(S > 0\), which is a Hilbert space with the following inner product and norm.
\begin{itemize}
\item[] $\bullet$ \  \(\mathcal{H}^S := {\rm L}^2(0, S; \mathbb{L}^2(\Omega)) \ni U(t), V(t)\),
\\[1mm]
\quad with inner product: \((U, V)_{\mathcal{H}^S} = \int_0^S (U, V)_{\mathbb{L}^2}^2 dt\) \ 
  and norm: \(||U||_{\mathcal{H}^S}^2 = (U, U)_{\mathcal{H}^S}\).
\end{itemize}

To apply the theory of parabolic equations, we write down each term of (CGL) in terms of subdifferential of some functional.

Let \({\rm H}\) be a real Hilbert space and denote by \(\Phi({\rm H})\) the set of all lower semi-continuous convex function \(\phi\) from \({\rm H}\) into \((-\infty, + \infty]\) such that the effective domain of \(\phi\) given by \({\rm D}(\phi) := \{u \in {\rm H} \mid \phi(u) < +\infty \}\) is not empty.
Then for \(\phi \in \Phi({\rm H})\), the subdifferential of \(\phi\) at \(u \in {\rm D}(\phi)\) is defined by
\begin{equation}
\label{subdiff}
\partial\phi(u) := \{f \in {\rm H} \mid (f, v - u)_{\rm H} \leq \phi(v) - \phi(u)\ \mbox{for all}\ v \in {\rm H}\}.
\end{equation}
Then \(\partial\phi\) becomes a (possibly multivalued) maximal monotone operator with domain \({\rm D}(\partial \phi) = \{u \in {\rm H} \mid \partial\phi(u) \neq \emptyset\}\). 
However for the arguments below, it suffices to deal with single-valued subdifferential operators.

Here we define two functionals on \(\mathbb{L}^2(\Omega)\).
\begin{align}
\label{fphi}
\varphi(U) &:=
\left\{
\begin{aligned}
&\frac{1}{2}\int_{\Omega}|\nabla U(x)|^2dx = \frac{1}{2}|\nabla U|_{\mathbb{L}^2}^2 &&\mbox{if}\ U \in \mathbb{H}_0^1(\Omega),
\\[1mm]
&+\infty &&\mbox{if}\ U \in \mathbb{L}^2(\Omega) \setminus \mathbb{H}_0^1(\Omega).%\mbox{else}.
\end{aligned}
\right.
\\[2mm]
\label{fpsir}
\psi_r(U) &:=
\left\{
\begin{aligned}
&\frac{1}{r}\int_{\Omega}|U(x)|^rdx = \frac{1}{r}|U|_{\mathbb{L}^r}^r
&&\mbox{if}\ U \in \mathbb{L}^r(\Omega) \quad (1< r< \infty ),
\\[1mm]
&+\infty &&\mbox{if}\ U \in \mathbb{L}^2(\Omega) \setminus \mathbb{L}^r(\Omega).%\text{else}.
\end{aligned}
\right.
\end{align}
Since these functionals are proper (\(\not\equiv +\infty\)), convex and lower semi-continuous, subdifferentials of these are given by
\begin{align}
\label{pphi}
\partial\varphi (U) &= -\Delta U, &{\rm D}(\partial\varphi) &= \left\{ U \in \mathbb{L}^2(\Omega) \mid U \in \mathbb{H}^2(\Omega) \cap \mathbb{H}_0^1(\Omega)\right\},
\\[2mm]
\label{ppsir}
\partial\psi_r (U) &= |U|^{r-2} U, &{\rm D}(\partial\psi_r) &= \left\{ U \in \mathbb{L}^2(\Omega) \mid \ |U|^{r - 2}U \in \mathbb{L}^2(\Omega)\right\}.
\end{align}

By the maximal monotonicity of subdifferential operators, we can consider their Yosida approximations.
Here we fix notations for resolvent operators and Yosida approximations, and collect their properties for later use.

Let \(\phi\) be a proper convex lower semi-continuous functional on a real Hilbert Space \({\rm H}\).
Since the subdifferential \(\partial\phi\) of \(\phi\) is maximal monotone in \({\rm H}\), we can define its resolvent \(J_\mu^\phi := (1 + \mu\partial\phi)^{-1}: {\rm H} \to {\rm D}(\partial\phi)\) for \(\mu > 0\) and the Yosida approximation of \(\partial\phi\) is given by \(\partial\phi_\mu := \partial\phi J_\mu^\phi\).
It is well known that the Yosida approximation of \(\partial\phi\) corresponds to the subdifferential of the Moreau-Yosida regularization \(\phi_\mu\) of \(\phi\), which is a Fr\'echet differentiable function given by
\begin{equation}\label{MYreg}
       \phi_\mu(u) = \inf_{v \in {\rm H}}\left\{\frac{1}{2\mu}|u - v|_{{\rm H}}^2 + \phi(v)\right\} 
                     = \frac{\mu}{2}|(\partial\phi)_{\mu}(u)|_{{\rm H}}^2 
                            + \phi(J_\mu^\phi(u)) \leq \phi(u),
\end{equation}
and the following inequality holds (see \cite{O1}, \cite{B1}, \cite{B2}):
\begin{equation}
\label{YA}
|\partial\phi_\mu(u)|_{\rm H} \leq |\partial\phi(u)|_{\rm H}\quad\mbox{for every}\ u \in {\rm D}(\partial\phi).
\end{equation}

Using these notations, we can rewrite our partial differential equation (CGL) in an evolution equation (ACGL) in \(\mathbb{L}^2(\Omega)\).
\begin{equation}
\label{ACGL-}
\tag*{(ACGL)}
     \frac{dU}{dt}(t) + (\lambda + \alpha I) ~\! \partial\varphi(U) 
         - (\kappa + \beta I) ~\! \partial\psi_q(U) - \gamma ~\! \partial\psi_2(U) = F(t),
\end{equation}
where the matrix \(I\) is given by
\begin{equation}\label{I}
   I =
       \begin{pmatrix}
           0 & 1 
\\
          -1 & 0
       \end{pmatrix}.
\end{equation}
This matrix \(I\) plays the same role in \(\mathbb{L}^2(\Omega)\) as the imaginary unit in terms of the inner product of $\mathbb{C}$, e.g.,
\begin{equation}
\label{I-i}
\begin{aligned}
      \Re ~\! (u, v)_{\mathbb{C}} = \Re ~\! u \bar{v} 
                                       & = (U, V)_{\mathbb{R}^2},
\\[2mm]
      \Im ~\! (u, v)_{\mathbb{C}} = \Im ~\! u \bar{v} 
                                       & = (U, IV)_{\mathbb{R}^2}, 
\end{aligned}
\end{equation}
where  \( (  \cdot, \cdot )_{\mathbb{C}}\) and \( (  \cdot, \cdot )_{\mathbb{R}^2}\) 
  denote the inner products in \( \mathbb{C} \) and \(\mathbb{R}^2\) respectively.
Here we collect basic properties of the matrix \(I\):
\begin{align}
\label{iu}
    & \mbox{imaginary unit} && : I^2 = -E,
\\[1mm]
\label{iso}
    & \mbox{isometry}       && : |U|_{\mathbb{L}^2} = |IU|_{\mathbb{L}^2},
\\[1mm]
\label{sksymm} 
    & \mbox{skew-symmetry}  && : (U, IV)_{\mathbb{L}^2} = -(IU, V)_{\mathbb{L}^2},
\\[1mm]
\label{comm}
    & \mbox{commutativity}  && : I \partial\varphi(U) 
                                    = \partial\varphi(IU),\quad I\partial\psi_r(U) 
                                      = \partial\psi_r(IU),
\\[1.5mm]
\label{orthoR}
    & \mbox{orthogonality in} \ \mathbb{R}^2  && : (U, IU)_{\mathbb{R}^2} = 0,  
\\[1.5mm]
\label{ortho}
    & \mbox{orthogonality in} \ \mathbb{L}^2  && : (U, IU)_{\mathbb{L}^2} 
                                                     = (U, I\partial\varphi(U))_{\mathbb{L}^2} 
                                                        = (U, I\partial\psi_r(U))_{\mathbb{L}^2} 
                                                          = 0,
\\[1.5mm]
\label{addortho1}
    & \mbox{orthogonality in} \ \mathbb{L}^2  && : (\partial\varphi_\mu(U), IU)_{\mathbb{L}^2} 
                                                     = 0 
                                            = (\partial\varphi_\mu(U), I\partial\varphi(U))_{\mathbb{L}^2},
\\[1.5mm]
\label{addortho2}
    & \mbox{orthogonality in} \ \mathbb{L}^2  && : (\partial\psi_{r, \mu}(U), IU)_{\mathbb{L}^2} 
                                                      = 0 
                                           = (\partial\psi_{r, \mu}(U), I\partial\psi_r(U))_{\mathbb{L}^2},
\\[1.5mm]
\label{bessel}
    & \mbox{Bessel's inequality}  && : (U, V)_{\mathbb{L}^2}^2 + (U, IV)_{\mathbb{L}^2}^2 
                                             \leq |U|_{\mathbb{L}^2}^2|V|_{\mathbb{L}^2}^2,
\end{align}
where \(E\) is the \(2 \times 2\) unit matrix,        
           \(\partial\varphi_\mu = (\partial\varphi)_\mu = \partial\varphi(1 + \mu \partial\varphi)^{-1}\)             and \(\partial\psi_{r, \mu} = (\partial\psi_r)_\mu 
                  = \partial\psi_r(1 + \mu \partial\psi_r)^{-1}\) 
                    denotes the Yosida approximations of \(\partial\varphi\) and \(\partial\psi_r\)                      respectively.

Though these properties can be proved by direct calculations, we give proofs of (\ref{addortho1}) and (\ref{addortho2}) for the sake of completeness.
%%%%%%%%%%%%%%%%%%%%%%%%%%%%%%%%%%%%%%%%%%%%%%%%%%
%%%%%%%%%% pf of addortho 1   %%%%%%%%%%%%%%%%%%%%
%%%%%%%%%%%%%%%%%%%%%%%%%%%%%%%%%%%%%%%%%%%%%%%%%%
\begin{Prf}[{\rm (\ref{addortho1})}]
Let \(V := (1 + \mu \partial\varphi)^{-1} U\).
Then by (\ref{sksymm}), (\ref{ortho}),  we get
\[
\begin{aligned}
   (\partial \varphi_\mu(U), IU)_{\mathbb{L}^2} 
       &= -(I\partial \varphi_\mu(U), U)_{\mathbb{L}^2} 
         = -(I \partial\varphi(V), (1+ \mu\partial\varphi)V)_{\mathbb{L}^2}\\ 
           &= -(I \partial\varphi(V), V)_{\mathbb{L}^2} = 0.
\end{aligned}
\]
 Furthermore, by virtue of (\ref{sksymm}), (\ref{ortho}) and the self-adjointness of 
$\partial \varphi(\cdot)$, we get
\[
\begin{aligned}
    &(\partial\varphi_\mu(U), I\partial\varphi(U))_{\mathbb{L}^2} 
       =  -(I\partial\varphi_\mu(U), \partial\varphi(U))_{\mathbb{L}^2} 
          = -(I\partial\varphi(V), \partial\varphi(U))_{\mathbb{L}^2}
\\
       =\mbox{}&  -\frac{1}{\mu}(I(U - V), \partial\varphi(U))_{\mathbb{L}^2}   
           = \frac{1}{\mu}(IV, \partial\varphi(U))_{\mathbb{L}^2}
             = \frac{1}{\mu}(I\partial\varphi(V), U)_{\mathbb{L}^2} 
\\
       =\mbox{}&  \frac{1}{\mu}(I\partial \varphi(V), V 
                 + \mu\partial\varphi(V))_{\mathbb{L}^2} 
           = \frac{1}{\mu}(I\partial\varphi(V), V)_{\mathbb{L}^2} 
             = 0.
\end{aligned}
\]
\end{Prf}
%%%%%%%%%%%%%%%%%%%%%%%%%%%%%%%%%%%%%%%%%%%%%%%%%%
%%%%%%%%%% end of pf of addortho 1   %%%%%%%%%%%%%
%%%%%%%%%%%%%%%%%%%%%%%%%%%%%%%%%%%%%%%%%%%%%%%%%%
%%%%%%%%%% pf of addortho 2   %%%%%%%%%%%%%%%%%%%%
%%%%%%%%%%%%%%%%%%%%%%%%%%%%%%%%%%%%%%%%%%%%%%%%%%
\begin{Prf}[{\rm (\ref{addortho2})}]
Let \(V := (1 + \mu\partial\psi_r)^{-1}U\).
By (\ref{sksymm}) and (\ref{ortho}), we obtain
\[
(\partial \psi_{r, \mu}(U), I U)_{\mathbb{L}^2} 
  = - ( I \partial \psi_{r, \mu}(U), U)_{\mathbb{L}^2} 
    = - ( I \partial\psi_r(V), V + \mu \partial\psi_r(V) )_{\mathbb{L}^2} 
      = 0.
\]
As for the second identity, we obtain by (\ref{sksymm}) and (\ref{orthoR}),
\[
\begin{aligned}
 & (\partial\psi_{r, \mu}(U), I\partial\psi_r(U))_{\mathbb{L}^2} 
  = - ( I\partial\psi_{r, \mu}(U), \partial\psi_r(U))_{\mathbb{L}^2}
\\[1mm]
 & = - \frac{1}{\mu}(I(U - V), \partial\psi_r(U))_{\mathbb{L}^2}
     = \frac{1}{\mu}(IV, \partial\psi_r(U))_{\mathbb{L}^2}
\\[1mm]
 & =  \frac{1}{\mu}\int_\Omega 
         (I V, |V + \mu\partial\psi_rV|_{\mathbb{R}^2}^{r - 2}
            (V + \mu\partial\psi_rV))_{\mathbb{R}^2} ~\! dx
\\[1mm]
 & =  \frac{1}{\mu}\int_\Omega
         (I V,  |V + \mu\partial\psi_rV|_{\mathbb{R}^2}^{r - 2}
            (V + \mu|V|_{\mathbb{R}^2}^{r - 2}V))_{\mathbb{R}^2} ~\! dx
              = 0,
\end{aligned}
\]
where we use temporal notation \(|\cdot|_{\mathbb{R}^2}\) for the length of vectors in \(\mathbb{R}^2\).
\end{Prf}
%%%%%%%%%%%%%%%%%%%%%%%%%%%%%%%%%%%%%%%%%%%%%%%%%%
%%%%%%%%%% end of pf of addortho 2   %%%%%%%%%%%%%
%%%%%%%%%%%%%%%%%%%%%%%%%%%%%%%%%%%%%%%%%%%%%%%%%%

Under these preparation we state our main results.
The first one is concerning the local well-posedness.
%%%%%%%%%%%%%%%%%%%%%%%%%%%%%%%%%%%%%%%%%%%%%%%%%%
%%%%%%%%%%%%%%%%%%%%%%%%%%%%%%%%%%%%%%%%%%%%%%%%%%
%%%%%%%   Theorem 1          %%%%%%%%%%%%%%%%%%%%%
%%%%%%%%%%%%%%%%%%%%%%%%%%%%%%%%%%%%%%%%%%%%%%%%%%
%%%%%%%%%%%%%%%%%%%%%%%%%%%%%%%%%%%%%%%%%%%%%%%%%%
\begin{Thm}[Local well-posedness in bounded domains]
\label{lwpndbd}
Let \(\Omega \subset \mathbb{R}^N\) be a bounded domain of \({\rm C}^2\)-regular class, \(F \in \mathcal{H}^T\) and \(2 < q < 2^*\) (subcritical), where
\[
   2^* =
     \begin{cases}
       + \infty & (N = 1, 2),
\\[2mm]
        \frac{2N}{N - 2} & (N \geq 3).
     \end{cases}
\]
Then for all \(U_0 \in \mathbb{H}_0^1(\Omega) = {\rm D}(\varphi)\), 
  there exist \( T_0 \in (0, T] \) and the unique function 
      \(U(t) \in {\rm C}([0, T_0]; \mathbb{L}^2(\Omega))\) satisfying:
\begin{enumerate}
\renewcommand{\labelenumi}{(\roman{enumi})}
\item \(U \in {\rm W}^{1, 2}(0, T_0; \mathbb{L}^2(\Omega))\),
\item \(U(t) \in {\rm D}(\partial\varphi) \subset {\rm D}(\partial\psi_q)\) for a.e. \(t \in [0, T_0]\) and satisfies (ACGL) for a.e. \(t \in [0, T_0]\),
\item \(\partial\varphi(U(\cdot)), \partial\psi_q(U(\cdot)) \in \mathcal{H}^{T_0}\).
\end{enumerate}
\end{Thm}
%%%%%%%%%%%%%%%%%%%%%%%%%%%%%%%%%%%%%%%%%%%%%%%%%%
Furthermore the following alternative on the maximal existence time of the solution holds:
%%%%%%%%%%%%%%%%%%%%%%%%%%%%%%%%%%%%%%%%%%%%%%%%%%
%%%%%%%%%    Theorem 2   %%%%%%%%%%%%%%%%%%%%%%%%%
%%%%%%%%%%%%%%%%%%%%%%%%%%%%%%%%%%%%%%%%%%%%%%%%%%
\begin{Thm}[Alternative]
\label{altbd}
Let \(T_m\) be the maximal existence time of the solution to (ACGL) 
  satisfying the regularity (i)-(iii) given in Theorem \ref{lwpndbd} 
for all $T_0 \in (0,T_m)$.
Then the following alternative on \(T_m\) holds:
\begin{itemize}
\item \(T_m = T\)\ \ \textbf{or}
\item \(T_m < T\) and \( \displaystyle\lim_{t \uparrow T_m}\varphi(U(t)) = +\infty\).
\end{itemize}
\end{Thm}
%%%%%%%%%%%%%%%%%%%%%%%%%%%%%%%%%%%%%%%%%%%%%%%%%%

In order to formulate the existence of small global solutions (see Theorem \ref{gebd}), we need to use the first eigenvalue \(\lambda_1 > 0 \) of \(-\Delta\) with the homogeneous Dirichlet boundary condition characterized by Poincar\'e's inequality:
\begin{equation}
\label{Poin}
        \psi_2(U) 
          \leq \lambda_1^{-1} \varphi(U),  \quad  \forall U \in \mathbb{H}_0^1(\Omega).
\end{equation}

For \(F \in {\rm L}^2(0, T; \mathbb{L}^2(\Omega))\), let \(\tilde{F}\) be the extention of \(F\) by zero to \((0, +\infty)\).
We set the following notation in order to scale the external force \(F\) in terms of \(\tilde{F}\)
\[
   \opnorm{F}^2 :=
                    \sup\left\{  \int_s^{s + 1}|\tilde{F}(t)|_{\mathbb{L}^2}^2 ~\! dt 
                                  \mid 0 \leq s < +\infty \right\}.
\]
%%%%%%%%%%%%%%%%%%%%%%%%%%%%%%%%%%%%%%%%%%%%%%%%%%
%%%%%%%%%%   Theorem 3    %%%%%%%%%%%%%%%%%%%%%%%%
%%%%%%%%%%%%%%%%%%%%%%%%%%%%%%%%%%%%%%%%%%%%%%%%%%
\begin{Thm}[Existence of small global solutions]
\label{gebd}
Let all the assumptions in Theorem \ref{lwpndbd} be satisfied and let 
 \(\gamma < \lambda ~\! \lambda_1\).
Then there exists a sufficiently small number \(r\) independent of \(T\) such that for all \(U_0 \in D(\varphi)\) and \(F \in {\rm L}^2(0, T; \mathbb{L}^2(\Omega))\) with \(\varphi(U_0) \leq r^2\) and \(\opnorm{F} \leq r\), every local solution given in Theorem \ref{lwpndbd} can be continued globally up to \([0, T]\).
\end{Thm}
%%%%%%%%%%%%%%%%%%%%%%%%%%%%%%%%%%%%%%%%%%%%%%%%%%
%%%%%%%%%%%%%%%%%%%%%%%%%%%%%%%%%%%%%%%%%%%%%%%%%%
%%%%%%%%%%   Aux eq   %%%%%%%%%%%%%%%%%%%%%%%%%%%%
%%%%%%%%%%%%%%%%%%%%%%%%%%%%%%%%%%%%%%%%%%%%%%%%%%
%%%%%%%%%%%%%%%%%%%%%%%%%%%%%%%%%%%%%%%%%%%%%%%%%%
\section{Solvability of Auxiliary Equation}
In this section, for any $S$ fixed in $(0,T]$, we consider the following auxiliary equation:
\begin{equation}
\tag*{(AE)\(^h\)}
\left\{
\begin{aligned}
    & \frac{dU^h}{dt}(t) + (\lambda + \alpha I) ~\! \partial \varphi(U^h) 
         - (\kappa + \beta I)  ~\! h(t) - \gamma ~\! U^h 
               = F(t)\\
&\hspace{80mm} t \in [0, S],
\\[2mm]
     & U^h(0) = U_0,
\end{aligned}
\right.
\end{equation}
which is (ACGL) with \(\partial\psi_q U(t)\) replaced by \(h \in \mathcal{H}^S\).

As for the global well-posedness for this auxiliary equation (AE)\(^h\), 
  we have:
\begin{Prop}
\label{gwpAEh}
Let \(\Omega \subset \mathbb{R}^N\) be a bounded domain of \({\rm C}^2\)-regular class, \(F \in \mathcal{H}^T\) and \(h \in \mathcal{H}^S\), \(0 < S \leq T\).
For all \(U_0 \in \mathbb{H}_0^1(\Omega) = D(\varphi)\), there exists the unique global solution \(U(t) \in {\rm C}([0, S]; \mathbb{L}^2(\Omega))\) satisfying:
\begin{enumerate}
\renewcommand{\labelenumi}{(\roman{enumi})}
\item \(U \in {\rm W}^{1, 2}(0, S; \mathbb{L}^2(\Omega))\),
\item \(U(t) \in {\rm D}(\partial\varphi)\) for a.e. \(t \in [0,S]\) and satisfies (AE)\(^h\) for a.e. \(t \in [0, S]\),
\item \(\partial\varphi(U(\cdot)) \in \mathcal{H}^S\).
\end{enumerate}
\end{Prop}
First we consider the following approximate equation:
\begin{equation}
\tag*{(AE)\(^h_\mu\)}
\left\{
   \begin{aligned}
   & \frac{dU_\mu}{dt}(t) + \lambda ~\! \partial\varphi(U_\mu) 
         + \alpha I \partial\varphi_\mu(U_\mu) - (\kappa + \beta I) ~\! h(t) - \gamma ~\!  U_\mu 
           = F(t), \\  
&\hspace{80mm}t \in [0, S],
\\[2mm]
   & U_\mu(0) = U_0,
   \end{aligned}
\right.
\end{equation}
which is (AE)\(^h\) with \(I\partial\varphi\) replaced by its Yosida approximation \(I\partial\varphi_\mu\).
By the standard theory of subdifferential operators, one can easily obtain 
the unique global solution for (AE)\(^h_\mu\) satisfying all properties (i)-(iii) 
stated in Proposition \ref{gwpAEh}, since the Yosida approximation 
\(\partial\varphi_\mu\) is Lipschitz continuous.

Here we are going to establish some a priori estimates for the solution \(U_\mu\) of 
(AE)\(^h_\mu\).
%%%%%%%%%%%%%%%%%%%%%%%%%%%%%%%%%%%%%%%%%%%%%%%%%%%%%%
%%%%%%%%%%   Lemma 1      %%%%%%%%%%%%%%%%%%%%%%%%%%%%
%%%%%%%%%%%%%%%%%%%%%%%%%%%%%%%%%%%%%%%%%%%%%%%%%%%%%%
\begin{Lem}[First Energy Estimate]
\label{1eAEm-}
Let \(U_\mu\) be the solution of (AE)\(^h_\mu\).
Then there exists \(C_1\) depending only on \(\lambda, \kappa, \beta, \gamma\), \(T\), \(|U_0|_{\mathbb{L}^2}, \|h\|_{\mathcal{H}^T}\) and \(\|F\|_{\mathcal{H}^T}\) such that
\begin{equation}
\label{1eAEm-1}
\sup_{t \in [0, S]}|U_\mu(t)|_{\mathbb{L}^2}^2 + \int_{0}^{S}\varphi(U_\mu(t)) ~\! dt \leq C_1.
\end{equation}
\end{Lem}
%%%%%%%%%%%%%%%%%%%%%%%%%%%%%%%%%%%%%%%%%%%%%%%%%%%%%%%
%%%%  Proof of lemma 1       %%%%%%%%%%%%%%%%%%%%%%%%%%
%%%%%%%%%%%%%%%%%%%%%%%%%%%%%%%%%%%%%%%%%%%%%%%%%%%%%%%
\begin{Prf}
Multiplying (AE)\(^h_\mu\) by \(U_\mu\) and by (\ref{addortho1}), we obtain
\begin{equation}
\label{AEm-Um}
   \begin{aligned}
     & \frac{1}{2}\frac{d}{dt}|U_\mu|_{\mathbb{L}^2}^2 + 2 \lambda ~\! \varphi(U_\mu)
\\[2mm]
      = \mbox{} & \gamma  ~\! |U_\mu|_{\mathbb{L}^2}^2 
          + \left((\kappa + \beta I)h + F, U_\mu\right)_{\mathbb{L}^2}
\\
     \leq\mbox{} & \left(\gamma_+ + \kappa^2 + \beta^2 
                     + \frac{1}{2}\right)|U_\mu|_{\mathbb{L}^2}^2 
                       + \frac{1}{2}|h|_{\mathbb{L}^2}^2 + \frac{1}{2}|F|_{\mathbb{L}^2}^2,
   \end{aligned}
\end{equation}
where we used the notation \(\gamma_+ := \max\{0, \gamma\}\) and Cauchy-Schwarz' inequality.
Integrating (\ref{AEm-Um}) on \((0, S)\), we obtain by Gronwall's inequality 
\[
   \begin{aligned}
     & \frac{1}{2}|U_\mu|_{\mathbb{L}^2}^2 + 2\lambda\int_0^t\varphi(U_\mu)d\tau
\\[1mm]
      \leq \mbox{} & \frac{1}{2}\left(|U_0|_{\mathbb{L}^2}^2 
                      + \|h\|_{\mathcal{H}^S}^2 + \|F\|_{\mathcal{H}^T}^2\right)
\\[1mm]
                   & + \left(2(\gamma_+ + \kappa^2 + \beta^2) + 1\right) 
                          \int_0^t\left(\frac{1}{2}|U_\mu|_{\mathbb{L}^2}^2 
                            + 2 \lambda\int_0^\tau\varphi(U_\mu)d\sigma\right)d\tau 
\\[1mm]
       \leq\mbox{} & \frac{1}{2}\left(|U_0|_{\mathbb{L}^2}^2 
                          + \|h\|_{\mathcal{H}^S}^2 
                            + \|F\|_{\mathcal{H}^T}^2 \right)
                               e^{\left(2(\gamma_+ + \kappa^2 + \beta^2) + 1\right)t}
\\[1mm]
       \leq\mbox{} & \frac{1}{2}\left(|U_0|_{\mathbb{L}^2}^2 
                          + \|h\|_{\mathcal{H}^S}^2 
                             + \|F\|_{\mathcal{H}^T}^2\right)
                               e^{\left(2(\gamma_+ + \kappa^2 + \beta^2) + 1\right)S},
   \end{aligned}
\]
which implies the desired estimate (\ref{1eAEm-1}).
\end{Prf}
%%%%%%%%%%%%%%%%%%%%%%%%%%%%%%%%%%%%%%%%%%%%%%%%%%%%%%
%%%%%%%    Lemma 2     %%%%%%%%%%%%%%%%%%%%%%%%%%%%%%%
%%%%%%%%%%%%%%%%%%%%%%%%%%%%%%%%%%%%%%%%%%%%%%%%%%%%%%
\begin{Lem}[Second Energy Estimates]
\label{2eAEm-}
Let \(U_\mu\) be the solution of (AE)\(^h_\mu\).
Then there exists \(C_2\) depending only on \(\lambda, \kappa, \beta, \gamma\), \(T\), \(|U_0|_{\mathbb{L}^2}, \varphi(U_0), \|h\|_{\mathcal{H}^S}\) and \(\|F\|_{\mathcal{H}^T}\) such that
\begin{equation}
\label{2eAEm-1}
\sup_{t \in [0, S]}\varphi(U_\mu(t)) + \int_{0}^{S}|\partial\varphi(U_\mu(t))|_{\mathbb{L}^2}^2dt + \int_{0}^{S}\left|\frac{dU_\mu}{dt}(t)\right|_{\mathbb{L}^2}^2dt \leq C_2.
\end{equation}
\end{Lem}
%%%%%%%%%%%%%%%%%%%%%%%%%%%%%%%%%%%%%%%%%%%%%%%%%%%%%%
%%%%%%%%%   Proof of lemma 2     %%%%%%%%%%%%%%%%%%%%%
%%%%%%%%%%%%%%%%%%%%%%%%%%%%%%%%%%%%%%%%%%%%%%%%%%%%%%
\begin{Prf}
Multiplying (AE)\(^h_\mu\) by \(\partial\varphi(U_\mu)\) and using (\ref{addortho1}) and Young's inequality, we obtain
\begin{equation}\label{ineq:varphiUmu}
\begin{aligned}
 & \frac{d}{dt}\varphi(U_\mu) 
      + \lambda ~\! |\partial\varphi(U_\mu)|_{\mathbb{L}^2}^2 
\\[2mm]
 & = 2 \gamma ~\!  \varphi(U_\mu) 
      + \left((\kappa + \beta I)h + F, \partial\varphi (U_\mu)\right)_{\mathbb{L}^2}
\\[1mm]
 & \leq 2 \gamma_+ ~\! \varphi(U_\mu) 
      + \frac{\kappa^2 + \beta^2}{\lambda}|h|_{\mathbb{L}^2}^2 
        + \frac{1}{\lambda}|F|_{\mathbb{L}^2}^2 
          + \frac{3\lambda}{4} ~\! |\partial\varphi(U_\mu)|_{\mathbb{L}^2}^2, 
\end{aligned}
\end{equation}
whence follows
\begin{equation}\label{AEm-pphiUm}
    \frac{d}{dt} \varphi(U_\mu) + \frac{\lambda}{4} ~\! |\partial\varphi(U_\mu)|_{\mathbb{L}^2}^2 
      \leq 2 \gamma_+ ~\! \varphi(U_\mu) + \frac{\kappa^2 + \beta^2}{\lambda}|h|_{\mathbb{L}^2}^2 
        + \frac{1}{\lambda}|F|_{\mathbb{L}^2}^2.
\end{equation}
Integrating (\ref{AEm-pphiUm}) on \((0, t)\) for \(t \in (0, S]\) and by Lemma \ref{1eAEm-}, we get
\begin{equation}
\label{IAEm-pphiUm}
\begin{aligned}
&\varphi(U_\mu) + \frac{\lambda}{4}\int_0^t|\partial\varphi(U_\mu(t))|_{\mathbb{L}^2}^2d\tau\\
&\leq \varphi(U_0) + 2\gamma_+\int_0^t\varphi(U_\mu)d\tau + \frac{\kappa^2 + \beta^2}{\lambda}\|h\|_{\mathcal{H}^S}^2 + \frac{1}{\lambda}\|F\|_{\mathcal{H}^T}^2\\
&\leq \varphi(U_0) + 2\gamma_+C_1 + \frac{\kappa^2 + \beta^2}{\lambda}\|h\|_{\mathcal{H}^S}^2 + \frac{1}{\lambda}\|F\|_{\mathcal{H}^T}^2\quad\ \mbox{for all}\ t \in (0, S].
\end{aligned}
\end{equation}
Thus from (\ref{IAEm-pphiUm}), (\ref{YA}) and (AE)\(^h_\mu\), we derive (\ref{2eAEm-1}).
\end{Prf}
%%%%%%%%%%%%%%%%%%%%%%%%%%%%%%%%%%%%%%%%%%%%%%%%%%%%%
%%%%   Proof of Proposition 4    %%%%%%%%%%%%%%%%%%%%
%%%%%%%%%%%%%%%%%%%%%%%%%%%%%%%%%%%%%%%%%%%%%%%%%%%%%
\begin{Prf}[Proposition \ref{gwpAEh}]
Let \(U_\mu\) be a solution of (AE)\(^h_\mu\).
First we show \(\{U_\mu\}_{\mu > 0}\) forms a Cauchy net in \({\rm C}([0, S]; \mathbb{L}^2(\Omega))\).
To this end, we multiply (AE)\(^h_\mu\)\(-\)(AE)\(_\nu\) by \(U_\mu - U_\nu\) to get
\begin{equation}\label{eqeq1}
   \begin{aligned}
    & \frac{1}{2}\frac{d}{dt}|U_\mu - U_\nu|_{\mathbb{L}^2}^2 
               + 2 \lambda ~\! \varphi(U_\mu - U_\nu)
\\[1mm]
    & = \gamma ~\! |U_\mu - U_\nu|_{\mathbb{L}^2}^2 
          - \alpha\left(I\partial\varphi_\mu U_\mu - I\partial\varphi_\nu U_\nu,
                                              U_\mu - U_\nu\right)_{\mathbb{L}^2}.
   \end{aligned}
\end{equation}
By the definition of Yosida approximation, it holds
\begin{equation}\label{eqeq2}
\begin{aligned}
&(I\partial\varphi_\mu U_\mu - I\partial\varphi_\nu U_\nu, U_\mu - U_\nu)_{\mathbb{L}^2}\\
%%&=\mbox{}(I\partial\varphi_\mu U_\mu - I\partial\varphi_\nu U_\nu, U_\mu - U_\nu)_{\mathbb{L}^2}\\
&=\mbox{}
\begin{aligned}[t]
&(I\partial\varphi_\mu U_\mu - I\partial\varphi_\nu U_\nu, (U_\mu - J_\mu U_\mu) - (U_\nu - J_\nu U_\nu))_{\mathbb{L}^2}\\
&+ (I(\partial\varphi_\mu U_\mu - \partial\varphi_\nu U_\nu), J_\mu U_\mu - J_\nu U_\nu)_{\mathbb{L}^2}
\end{aligned}\\
&=\mbox{}
\begin{aligned}[t]
&(I\partial\varphi_\mu U_\mu - I\partial\varphi_\nu U_\nu, \mu\partial_\mu U_\mu - \nu\partial\varphi_\nu U_\nu)_{\mathbb{L}^2}\\
&+(I(\partial\varphi J_\mu U_\mu - \partial\varphi J_\nu U_\nu), J_\mu U_\mu - J_\nu U_\nu)_{\mathbb{L}^2}.
\end{aligned}
\end{aligned}
\end{equation}
We make use of the linearity of \(\partial\varphi\) and \eqref{ortho} to obtain
\begin{equation}\label{eqeq3}
\begin{aligned}
&(I(\partial\varphi J_\mu U_\mu - \partial\varphi J_\nu U_\nu), J_\mu U_\mu - J_\nu U_\nu)_{\mathbb{L}^2}\\
&=\mbox{}(I\partial\varphi(J_\mu U_\mu - J_\nu U_\nu), J_\mu U_\mu - J_\nu U_\nu)_{\mathbb{L}^2} = 0.
\end{aligned}
\end{equation}
%%
%Applying K\=omura's trick, we obtain
%%
Combining \eqref{eqeq1}-\eqref{eqeq3}, we have
\[
   \begin{aligned}
     & \frac{1}{2}\frac{d}{dt}|U_\mu - U_\nu|_{\mathbb{L}^2}^2
\\[1mm]
      \leq \mbox{} & \gamma_+|U_\nu - U_\mu|_{\mathbb{L}^2}^2 
           + |\alpha| \left\{ \mu ~\! |\partial\varphi_\nu(U_\nu)|_{\mathbb{L}^2}|                    
              \partial\varphi_\mu(U_\mu)|_{\mathbb{L}^2} 
                + \nu ~\! |\partial\varphi_\mu(U_\mu)|_{\mathbb{L}^2}
                     |\partial\varphi_\nu(U_\nu)|_{\mathbb{L}^2} \right\}
\\[1mm]
     \leq \mbox{} & \gamma_+|U_\nu - U_\mu|_{\mathbb{L}^2}^2 
                    + \frac{|\alpha|}{2}(\mu + \nu) 
                        \left\{|\partial\varphi(U_\mu)|_{\mathbb{L}^2}^2 
                            + |\partial \varphi(U_\nu)|_{\mathbb{L}^2}^2\right\}.
   \end{aligned}
\]
Thus Gronwall's inequality yields
\[
    |U_\mu(t) - U_\nu(t)|_{\mathbb{L}^2}^2 
        \leq |\alpha|(\mu + \nu) ~\! e^{2\gamma_+t} 
               \int_0^t\left\{|\partial\varphi(U_\mu(s))|_{\mathbb{L}^2}^2 
                 + |\partial \varphi(U_\nu(s))|_{\mathbb{L}^2}^2\right\}ds,
\]
for all \(t \in [0, S]\).
Then by Lemma \ref{2eAEm-}, we have
\[
\sup_{t \in [0, T]} |U_\mu(t) - U_\nu(t)|_{\mathbb{L}^2} \leq e^{\gamma_+T}\sqrt{2 C_2 |\alpha|(\mu + \nu)},
\]
which assures that \(\{ U_\mu \}_{\mu > 0}\) forms a Cauchy net in \({\rm C}([0, S]; \mathbb{L}^2(\Omega))\). 

   Now let \(U_{\mu} \rightarrow U\) in \({\rm C}([0, S]; \mathbb{L}^2(\Omega))\) as \(\mu \to 0\).
By Lemma \ref{2eAEm-}, \(\{\frac{d}{dt}U_{\mu}\}_{\mu > 0}\) and \(\{\partial\varphi(U_{\mu})\}_{\mu > 0}\) are bounded in \({\rm L}^2(0, S; \mathbb{L}^2(\Omega))\).
Hence by the demiclosedness of \(\frac{d}{dt}\) and \(\partial\varphi\) we have
\[
\begin{aligned}
\frac{dU_{\nu_n}}{dt} &\rightharpoonup \frac{dU}{dt} &\mbox{weakly in}\ {\rm L}^2(0, T; \mathbb{L}^2(\Omega)),\\
\partial\varphi(U_{\nu_n}) &\rightharpoonup \partial\varphi(U) &\mbox{weakly in}\ {\rm L}^2(0, T; \mathbb{L}^2(\Omega)),
\end{aligned}
\]
for some sequence \(\{\nu_n\}_{n \in \mathbb{N}}\) such that \(\nu_n \to 0\) as \(n \to \infty\).
We can also find
\[
J_{\nu_n}^\varphi U_{\nu_n} \rightarrow U\ \mbox{strongly in}\ {\rm L}^2(0, T; \mathbb{L}^2(\Omega)).
\]
Indeed by the definition of Yosida approximation, it holds that
\begin{equation}
\label{resolvent}
\begin{aligned}
||U_{\nu_n} - J_{\nu_n}^\varphi U_{\nu_n}||_{\mathcal{H}^T}^2 &= \int_0^T |U_{\nu_n}(s) - J_{\nu_n}^\varphi U_{\nu_n}(s)|_{\mathbb{L}^2}^2ds\\
&= \nu_n^2 \int_0^T |\partial\varphi_{\nu_n}(U_{\nu_n}(s))|_{\mathbb{L}^2}^2ds \leq C_2 \nu_n^2 \rightarrow 0\quad\mbox{as}\ n \rightarrow \infty.
\end{aligned}
\end{equation}
Then since \(\partial\varphi_\nu(U_\nu) = \partial\varphi(J_\nu^\varphi U_\nu)\), by the demiclosedness of \(\partial\varphi\) we find that \(U\) satisfies 
\[
   \frac{dU}{dt} + (\lambda + \alpha I) ~\! \partial \varphi(U) 
           - (\kappa + \beta I) ~\! h(t) - \gamma ~\! U 
              = F \quad \mbox{in} \ {\rm L}^2(0, T; \mathbb{L}^2(\Omega)),
\]
i.e., \(U\) is the desired solution of (AE)\(^h\).
\end{Prf}
%%%%%%%%%%%%%%%%%%%%%%%%%%%%%%%%%%%%%%%%%%%%%%%%%%
%%%%%%%%%%%%%%%%%%%%%%%%%%%%%%%%%%%%%%%%%%%%%%%%%%
%%%%%%%%%%   Pf of Thm 1 (ex)   %%%%%%%%%%%%%%%%%%
%%%%%%%%%%%%%%%%%%%%%%%%%%%%%%%%%%%%%%%%%%%%%%%%%%
%%%%%%%%%%%%%%%%%%%%%%%%%%%%%%%%%%%%%%%%%%%%%%%%%%
\section{Proof of Theorem \ref{lwpndbd} (Existence)}
Before proceeding to the proof of Theorem \ref{lwpndbd}, we establish some a priori estimates for the unique solutions \(U^h\) of auxiliary equations (AE)\(^h\), which are given in Proposition \ref{gwpAEh}.
First fix a constant \(R\) such as
\begin{equation}
\label{myR}
R := \max\left\{\frac{1}{2}|U_0|_{\mathbb{L}^2}^2 + \varphi(U_0) + \frac{1}{\lambda}\|F\|_{\mathcal{H}^T}^2, 1\right\}
\end{equation}
and define the closed convex subset $K_R^S$ of $\mathcal{H}^S$ by
\begin{equation}
\label{def:KRS}
 K_R^S := \{ ~\! h \in {\cal H}^S ~\! ; ~\! \|h\|_{\mathcal{H}^S} \leq R ~\! \}.
\end{equation}
%%%%%%%%%%%%%%%%%%%%%%%%%%%%%%%%%%%%%%%%%%%%%%%%%%
%%%%%%%%    Lemma 3     %%%%%%%%%%%%%%%%%%%%%%%%%%
%%%%%%%%%%%%%%%%%%%%%%%%%%%%%%%%%%%%%%%%%%%%%%%%%%
\begin{Lem}[First Energy Estimate]
\label{1eAEh}
Let $h \in K_R^S$ and \(U^h\) be the unique solution of (AE)\(^h\).
Then there exists \(C_1\) depending only on \(\lambda, \kappa, \beta\) and \(\gamma\) such that
\begin{equation}\label{1eAEh1}
     \sup_{t \in [0, S]}|U^h(t)|_{\mathbb{L}^2}^2 + \int_{0}^{S}\varphi(U^h(t)) ~\! dt 
       \leq C_1 R.
\end{equation}
\end{Lem}
%%%%%%%%%%%%%%%%%%%%%%%%%%%%%%%%%%%%%%%%%%%%%%%%%%%
%%%%%%%%%   Proof of Lemma 3     %%%%%%%%%%%%%%%%%%
%%%%%%%%%%%%%%%%%%%%%%%%%%%%%%%%%%%%%%%%%%%%%%%%%%%
\begin{Prf}
Multiply (AE)\(^h\) by \(U^h\), then by \eqref{ortho} we get ( see \eqref{AEm-Um})
\begin{equation}\label{AEhUh}
   \begin{aligned}
     & \frac{1}{2}\frac{d}{dt}|U^h|_{\mathbb{L}^2}^2 + 2 \lambda ~\! \varphi(U^h)
\\[2mm]
    = \mbox{} & \gamma ~\! |U^h|_{\mathbb{L}^2}^2 
                + \left((\kappa + \beta I)h + F, U^h\right)_{\mathbb{L}^2}
\\[2mm]
    \leq \mbox{} & \frac{4\gamma_+ + \kappa^2 + \beta^2 
                     + \lambda}{4} ~\! |U^h|_{\mathbb{L}^2}^2 
                       + |h|_{\mathbb{L}^2}^2 + \frac{1}{\lambda}|F|_{\mathbb{L}^2}^2,
   \end{aligned}
\end{equation}
where we used the notation \(\gamma_+ := \max\{0, \gamma\}\) and the Cauchy-Schwarz inequality.

Integrating (\ref{AEhUh}) on \((0, S)\) and noting the fact that 
$h \in K_R^S$, we obtain
\begin{equation}
\label{AEhUh1}
\begin{aligned}
&\frac{1}{2}|U^h(t)|_{\mathbb{L}^2}^2 + 2\lambda\int_0^t\varphi(U^h)d\tau\\
\leq\mbox{}& \frac{1}{2}|U_0|_{\mathbb{L}^2}^2 + \frac{1}{\lambda}\|F\|_{\mathcal{H}^T}^2 + \|h\|_{\mathcal{H}^S}^2\\
&+ \frac{4\gamma_+ + \kappa^2 + \beta^2 + \lambda}{2}\int_0^t\left(\frac{1}{2}|U^h|_{\mathbb{L}^2}^2 + 2\lambda\int_0^\tau\varphi(U^h)d\sigma\right)d\tau\\
\leq\mbox{}& 2R + \frac{4\gamma_+ + \kappa^2 + \beta^2 + \lambda}{2}\int_0^t\left(\frac{1}{2}|U^h|_{\mathbb{L}^2}^2 + 2\lambda\int_0^\tau\varphi(U^h)d\sigma\right)d\tau.
\end{aligned}
\end{equation}
We apply Gronwall's inequality to (\ref{AEhUh1}) to get
\[
\begin{aligned}
\frac{1}{2}|U^h(t)|_{\mathbb{L}^2}^2 + 2\lambda\int_0^t\varphi(U^h)d\tau &\leq 2Re^{\frac{4\gamma_+ + \kappa^2 + \beta^2 + \lambda}{2}t}\\
&\leq 2e^{\frac{4\gamma_+ + \kappa^2 + \beta^2 + \lambda}{2}S}R\quad\mbox{for all}\ t \in [0, S],
\end{aligned}
\]
which implies the desired estimate (\ref{1eAEh1}).
\end{Prf}
%%%%%%%%%%%%%%%%%%%%%%%%%%%%%%%%%%%%%%%%%%%%%%%%%%%%%%%%%%
%%%%%%%%  Lemma 4     %%%%%%%%%%%%%%%%%%%%%%%%%%%%%%%%%%%%
%%%%%%%%%%%%%%%%%%%%%%%%%%%%%%%%%%%%%%%%%%%%%%%%%%%%%%%%%%
\begin{Lem}[Second Energy Estimates]
\label{2eAEh}
Let \(U^h\) be the solution of (AE)\(^h\).
Then there exists \(C_2\) depending only on \(\lambda, \kappa, \beta\) and \(\gamma\) such that
\begin{equation}
\label{2eAEh1}
\sup_{t \in [0, S]}\varphi(U^h(t)) + \int_{0}^{S}|\partial\varphi(U^h(t))|_{\mathbb{L}^2}^2dt + \int_{0}^{S}\left|\frac{dU^h}{dt}(t)\right|_{\mathbb{L}^2}^2dt \leq C_2R.
\end{equation}
\end{Lem}
%%%%%%%%%%%%%%%%%%%%%%%%%%%%%%%%%%%%%%%%%%%%%%%%%%%%%%%%%
%%%%%% Proof of Lemma 4     %%%%%%%%%%%%%%%%%%%%%%%%%%%%%
%%%%%%%%%%%%%%%%%%%%%%%%%%%%%%%%%%%%%%%%%%%%%%%%%%%%%%%%%
\begin{Prf}
Multiplying (AE)\(^h\) by \(\partial\varphi(U^h)\) and using (\ref{addortho1}), we obtain (see \eqref{ineq:varphiUmu} ) 
\[
   \begin{aligned}
     & \frac{d}{dt} \varphi(U_\mu) + \lambda ~\! |\partial\varphi(U_\mu(t))|_{\mathbb{L}^2}^2
\\[1mm]
     & = 2 \gamma ~\! \varphi(U^h) + \left((\kappa + \beta I)h + F, 
                    \partial\varphi(U^h)\right)_{\mathbb{L}^2}
\\[1mm]
     & \leq 2 \gamma_+ ~\! \varphi(U^h) + \frac{\kappa^2 + \beta^2}{\lambda}|h|_{\mathbb{L}^2}^2 
               + \frac{1}{\lambda}|F|_{\mathbb{L}^2}^2 
                  + \frac{3\lambda}{4} ~\! |\partial\varphi(U^h(t))|_{\mathbb{L}^2}^2,
   \end{aligned}
\]
whence follows
\begin{equation}\label{AEhpphiUh}
    \frac{d}{dt} \varphi(U^h) + \frac{\lambda}{4} ~\! |\partial\varphi(U^h(t))|_{\mathbb{L}^2}^2 
       \leq 2 \gamma_+ ~\! \varphi(U^h) + \frac{\kappa^2 + \beta^2}{\lambda}|h|_{\mathbb{L}^2}^2 
                + \frac{1}{\lambda}|F|_{\mathbb{L}^2}^2.
\end{equation}
Integrating (\ref{AEhpphiUh}) on \((0, t)\) for \(t \in (0, S]\) and by Lemma \ref{1eAEh}, we get
\begin{equation}\label{IAEhpphiUh}
   \begin{aligned}
     &  \varphi(U^h) + \frac{\lambda}{4} \int_0^t\!|\partial\varphi(U^h(t))|_{\mathbb{L}^2}^2 ~\! d\tau
\\[1mm]
     &  \leq \varphi(U_0) + \frac{1}{\lambda}\|F\|_{\mathcal{H}^T}^2 
               + \frac{\kappa^2 + \beta^2}{\lambda}\|h\|_{\mathcal{H}^S}^2 
                  + 2\gamma_+\int_0^t\varphi(U^h)d ~\! \tau
\\[1mm]
     &  \leq \left(1 + \frac{\kappa^2 + \beta^2}{\lambda} + 2 \gamma_+C_1 \right) R 
               \quad  \ \mbox{for all}\ t \in (0, S].
   \end{aligned}
\end{equation}
Thus from (\ref{IAEhpphiUh}) and (AE)\(^h\), we derive (\ref{2eAEh1}).
\end{Prf}
%%%%%%%%%%%%%%%%%%%%%%%%%%%%%%%%%%%%%%%%%%%%%%%%%
Now we are ready to prove the existence part of Theorem \ref{lwpndbd}.

%%%%%%%%%%%%%%%%%%%%%%%%%%%%%%%%%%%%%%%%%%%%%%%%%%%%%
%%%%% Proof of Theorem 1     %%%%%%%%%%%%%%%%%%%%%%%%
%%%%%%%%%%%%%%%%%%%%%%%%%%%%%%%%%%%%%%%%%%%%%%%%%%%%%
\begin{Prf}[Theorem \ref{lwpndbd} (Existence)]
Let $K_R^S$ be the closed convex subset of $\mathcal{H}^S$ defined 
by \eqref{def:KRS} and we introduce a mapping $\mathcal{F}$ by the following correspondence:
\begin{equation}
\label{mathcalF}
\mathcal{F}: \mathcal{H}^S \ni h(t) \mapsto \mathcal{F}(h(t)) := \partial\psi_q(U^h) \in \mathcal{H}^S,
\end{equation}
where \(U^h\) is the unique solution of (AE)\(^h\).

First we show that \(\mathcal{F}\) maps \(K_R^S\) into itself for a sufficiently small $S \in (0,T]$.
By the Gagliardo{-}Nirenberg{-}Sobolev inequality, for any $q \in (2,2^*)$ there 
exists a constant \(C_{GN}\) such that 
\begin{equation}
\label{GNSH2}
   |\partial\psi_q(U^h)|_{\mathbb{L}^2}^2 
    = |U^h|_{\mathbb{L}^{2(q-1)}}^{2(q-1)} 
       \leq C_{\rm GN} ~\! |U^h|_{\mathbb{H}^2}^{2(1 - \xi)(q - 1)} 
               |U^h|_{\mathbb{L}^{2^*}}^{2\xi(q - 1)} 
                       \quad \forall U \in \mathbb{H}^2(\Omega), 
\end{equation}
where parameter \(\xi\) satisfies
\[
\frac{1}{2(q - 1)} = \left(\frac{1}{2} - \frac{2}{N}\right)(1 - \xi) + \left(\frac{1}{2} - \frac{1}{N}\right)\xi.
\]
We apply the elliptic estimate to (\ref{GNSH2}) to obtain
\begin{equation}\label{GNSE}
   |U^h|_{\mathbb{H}^2}^{2(1 - \xi)(q - 1)}
        |U^h|_{\mathbb{L}^{2^*}}^{2\xi(q - 1)}   
            \leq C ~\! |\partial\varphi(U^h)|_{\mathbb{L}^2}^{2(1 - \xi)(q - 1)}      
                ~\! \varphi(U^h)^{\xi(q - 1)},
\end{equation}
where \(C\) denotes some embedding constant.
Our assumption on \(q\) being Sobolev subcritical assures \((1 - \xi)(q - 1) < 1\).
Then by Young's inequality, for arbitrary \(\varepsilon > 0\), 
there exists a constant $C_\varepsilon$ such that 
\begin{equation}\label{GNSEY}
  \begin{aligned}
   &|\partial\varphi(U^h)|_{\mathbb{L}^2}^{2(1 - \xi)(q - 1)} 
         \varphi(U^h)^{\xi(q - 1)}    
           \leq \varepsilon ~\! |\partial\varphi(U^h)|_{\mathbb{L}^2}^2
                    + C_\varepsilon\varphi(U^h)^\chi, 
  \\[1mm]
    & \chi = \frac{\xi (q-1)}{1 - (q-1)(1-\xi)} >1.
  \end{aligned} 
\end{equation}
Here we note that $\chi>1$ if and only if $q>2$. 
   Hence by (\ref{GNSH2}), (\ref{GNSE}) and (\ref{GNSEY}), we get
\begin{equation}
\label{ppsiqUlL}
|\partial\psi_q(U^h)|_{\mathbb{L}^2}^2 \leq \varepsilon\left(|\partial\varphi(U^h)|_{\mathbb{L}^2}^2 + |U^h|_{\mathbb{L}^2}^2\right) + C_\varepsilon\varphi(U^h)^\chi.
\end{equation}
Integration of (\ref{ppsiqUlL}) on \([0, S]\) together with (\ref{2eAEh1}) gives
\[
   \begin{aligned}
      \int_{0}^{S}|\partial\psi_q(U^h)|_{\mathbb{L}^2}^2 dt 
       & \leq \varepsilon\int_{0}^{S}\left(|\partial\varphi(U^h)|_{\mathbb{L}^2}^2 
                + |U^h|_{\mathbb{L}^2}^2\right)dt + C_\varepsilon \int_{0}^{S}\!\varphi(U^h)^{\chi_2}dt
\\[2mm]
       & \leq \varepsilon ~\! C_2 ~\! R + M_\varepsilon(R) ~\! S,
\end{aligned}
\]
where \(M_\varepsilon(\cdot)\) denotes a non-decreasing function depending on \(\varepsilon\).

First fix \(\varepsilon := \frac{1}{2C_2}\) and then define \(S\) by
\begin{equation}
\label{myS}
S := \min\left\{T, \frac{R}{2M_\varepsilon(R)}\right\}.
\end{equation}
Then \(\int_{0}^{S}|\partial\psi_q(U^h)|_{\mathbb{L}^2}^2dt = \int_{0}^{S}|\mathcal{F}(h)|_{\mathbb{L}^2}^2dt \leq R\), that is \(\mathcal{F}\) maps \(K_R^S\) into itself.

Next we prove the weak continuity of \(\mathcal{F}\).
% continuity is a local property, we could focus on bounded neighbourhoods, which are metrizable because %\(\mathcal{H}^S = {\rm L}^2(0, S; \mathbb{L}^2(\Omega))\) is a separable Hilbert space.
Let \(\{h_n\}_{n \in \mathbb{N}}\) be a sequence in \(\mathcal{H}^S\) such that
\[
h_n \rightharpoonup h\ \mbox{weakly in}\ {\rm L}^2(0, S; \mathbb{L}^2(\Omega)),
\]
and \(U^{h_n}, U^h\) be unique solutions of (AE)\(^{h_n}\) and (AE)\(^h\) respectively.
Lemma \ref{2eAEh} assures the equi-continuity of \(\left\{U^{h_n}(t)\right\}_{n \in \mathbb{N}}\). In fact we have by (\ref{2eAEh1})
\[
   \begin{aligned}
    |U^{h_n}(t) - U^{h_n}(s)|_{\mathbb{L}^2} 
       = \left| \int_{s}^{t}\frac{dU^{h_n}}{d\tau}(\tau)d ~\! \tau \right|_{\mathbb{L}^2}
         & \leq \int_{t}^{s} \left| \frac{dU^{h_n}}{d\tau}(\tau) \right|_{\mathbb{L}^2} \!\! d\tau
\\
         & \leq \left(\int_{s}^{t} \left|\frac{dU^{h_n}}{d\tau}(\tau) \right|_{\mathbb{L}^2}^2 d \tau 
                        \right)^{\frac{1}{2}} \left(\int_{s}^{t} 1 ~\! d\tau \right)^{\frac{1}{2}}
\\
         & \leq \sqrt{C_2 R} \sqrt{t - s}.
   \end{aligned}
\]
Lemma \ref{2eAEh} and Rellich{-}Kondrachov theorem assure that the set \(\left\{U^{h_n}(t)\right\}_{n \in \mathbb{N}}\) is relatively compact in \(\mathbb{L}^2(\Omega)\) for all \(t \in [0, S]\).
 By Ascoli's Theorem and Lemma \ref{2eAEh}, there exists a subsequence \(\{h_{n'}\}_{n' \in \mathbb{N}} \subset \{h_n\}_{n \in \mathbb{N}}\) and \(U \in {\rm C}([0, S]; \mathbb{L}^2(\Omega))\) such that
\begin{align}
\label{Uh}
 U_{h_{n'}} &\rightarrow U
   &&%\begin{aligned}
%    &
\mbox{strongly in}\ {\rm C}(0, T; \mathbb{L}^2(\Omega))%\\
%     &\mbox{and strongly in}\ {\rm L}^2(0,T;\mathbb{L}^2(\Omega)),
 %    \end{aligned}
%\\[1mm]
\\[2mm]
\label{D-tUh}
\frac{dU^{h_{n'}}}{dt} &\rightharpoonup \frac{dU}{dt}&&\mbox{weakly\ in}\ {\rm L}^2(0,T;\mathbb{L}^2(\Omega)),
\\[2mm]
\label{delphiUh}
\partial\varphi(U^{h_{n'}}) &\rightharpoonup \partial\varphi(U)&&\mbox{weakly in}\ {\rm L}^2(0, T; \mathbb{L}^2(\Omega)),
\\[2mm]
\label{delpsiUh}
\partial\psi_q(U^{h_{n'}})&\rightharpoonup \partial\psi_q(U)&&\mbox{weakly in}\ {\rm L}^2(0, T; \mathbb{L}^2(\Omega)),
\end{align}
here we used the demiclosedness of \(\frac{d}{dt}\), \(\partial \varphi\) and \(\partial\psi_q\) in \({\rm L}^2(0, T; \mathbb{L}^2(\Omega))\) in (\ref{D-tUh}), (\ref{delphiUh}) and (\ref{delpsiUh}).
%%%%%%%
Thus \(U\) satisfies the following equation:
\[
    \frac{dU}{dt} + (\lambda + \alpha I) ~\! \partial\varphi(U) 
                   - (\kappa + \beta I) ~\! h - \gamma ~\!  U = F,
\]
i.e., \(U\) coincides with its unique solution \(U^h\).
%%%%%%
 Since the above argument does not depend on the choice subsequences, 
 we conclude that 
\[\mathcal{F}(h_n) = \partial\psi_q(U^{h_n}) \rightharpoonup \partial\psi_q(U) = \partial\psi_q(U^h) = \mathcal{F}(h),
\]
whence follows the weak continuity of \(\mathcal{F}\).

Now, we can apply Schauder-Tychonoff's fixed point theorem on \(\mathcal{F}\) and \(K_R^S\) to obtain a fixed point \(h\), i.e., \(h\) satisfies
\begin{equation}
\label{FP}
h = \mathcal{F}(h) = \partial\psi_q(U^h).
\end{equation}
By (\ref{FP}), the corresponding solution \(U^h\) satisfies:
\begin{equation}\label{concl}
   \begin{aligned}
    & \frac{dU^h}{dt} + (\lambda + \alpha I)~\! \partial\varphi(U^h) 
                 - (\kappa + \beta I) ~\! h - \gamma ~\! U^h
\\[1mm]
    = \mbox{} & \frac{dU^h}{dt} + (\lambda + \alpha I) ~\! \partial\varphi(U^h) 
              - (\kappa + \beta I) ~\! \partial\psi_q(U^h) - \gamma ~\! U^h = F,
\end{aligned}
\end{equation}
which means \(U^h\) is the desired solution of (ACGL).
\end{Prf}

%%%%%%%%%%%%%%%%%%%%%%%%%%%%%%%%%%%%%%%%%%%%%%%%%%
%%%%%%%%%%%%%%%%%%%%%%%%%%%%%%%%%%%%%%%%%%%%%%%%%%
%%%%%%%%%%   Pf of Thm 2   %%%%%%%%%%%%%%%%%%%%%%%
%%%%%%%%%%%%%%%%%%%%%%%%%%%%%%%%%%%%%%%%%%%%%%%%%%
%%%%%%%%%%%%%%%%%%%%%%%%%%%%%%%%%%%%%%%%%%%%%%%%%%
\section{Proof of Theorem \ref{altbd}}
Before showing the uniqueness of the solution for (ACGL), we prove Theorem \ref{altbd}.

Let \(T_m\) be the maximal existence time of a solution of (ACGL), i.e., 
\begin{align*}
 T_m := \sup  ~\! \{ ~\! S > 0 \  \mid \ & \exists\ \mbox{a solution of (ACGL) on}\ [0, S]
\\
         &    \text{satisfying (i)-(iii) of Theorem 1 with} \ T_0 = S ~\! \}.
\end{align*}
%%%%%%%%%%%%%%%%%%%%%%%%%%%%%%%%%%%%%%%%%%%%%%%%%%
%%%% Proof of Theorem 2    %%%%%%%%%%%%%%%%%%%%%%%
%%%%%%%%%%%%%%%%%%%%%%%%%%%%%%%%%%%%%%%%%%%%%%%%%%
\begin{Prf}[Theorem \ref{altbd}]
We rely on proof by contradiction.
Assume \(T_m < T\) and the assertion \(\lim_{t \uparrow T_m}\varphi(U(t)) = +\infty\) does not hold.
Then there exists monotonically increasing sequence \(t_n \uparrow T_m\) such that \(\varphi(U(t_n)) \leq C\) holds for all \(n \in \mathbb{N}\). 
We repeat the same argument as before with \(U(0)\) replaced by \(U(t_n)\) to assure the existence of \(\sigma > 0\) independent of \(n\) such that a solution of (ACGL) exists on \([t_n, t_n + \sigma]\).
Recalling the definition (\ref{myR}) of \(R\), we define
\[
\begin{aligned}
\rho := \max\left\{C(\lambda_1^{-1} + 1) + \frac{1}{\lambda}\|F\|_{\mathcal{H}}^2, 1\right\} \geq 1.
\end{aligned}
\]
Then by Poincar\'e's inequality, it holds for all \(n \in \mathbb{N}\)
\[
\rho \geq \frac{1}{\lambda}\|F\|_{\mathcal{H}}^2+\textstyle{\frac{1}{2}}|U(t_n)|_{\mathbb{L}^2}^2 + \varphi(U(t_n)).
\]
Additionally we define \(\sigma\) by ( see \eqref{myS} )
\[
\sigma := \min\left\{T - T_m, \frac{\rho}{2M_\varepsilon(\rho)}\right\},
\]
which is independent of \(n\).
We can deduce \(\mathcal{F}\) maps \(K_\rho^\sigma\) into itself in the same way as before.
Thus we can construct solution on \([t_n, t_n + \sigma]\) applying Schauder-Tychonoff's fixed point theorem again.

Since \(\{t_n\}_{n \in \mathbb{N}}\) converges to \(T_m\), there exists \(N_0 \in \mathbb{N}\) such that for all \(n \geq N_0\), it holds that \(T_m < t_n + \frac{\sigma}{2}\).
This means the local solution can be extended up to \(\left[0, T_m + \frac{\sigma}{2}\right]\), 
which contradicts the definition of \(T_m\).
\end{Prf}

%%%%%%%%%%%%%%%%%%%%%%%%%%%%%%%%%%%%%%%%%%%%%%%%%%
%%%%%%%%%%%%%%%%%%%%%%%%%%%%%%%%%%%%%%%%%%%%%%%%%%
%%%%%%%%%%   Pf of Thm 1 (Uniqueness)   %%%%%%%%%%
%%%%%%%%%%%%%%%%%%%%%%%%%%%%%%%%%%%%%%%%%%%%%%%%%%
%%%%%%%%%%%%%%%%%%%%%%%%%%%%%%%%%%%%%%%%%%%%%%%%%%
\section{Proof of Theorem \ref{lwpndbd} (Uniqueness)}
Here we give a proof for the uniqueness of solution of (ACGL). 
 We first prepare the following lemma.
%%%%%%%%%%%%%%%%%%%%%%%%%%%%%%%%%%%%%%%%%%%%%%%%%
%%%%%%     Lemma 5      %%%%%%%%%%%%%%%%%%%%%%%%%
%%%%%%%%%%%%%%%%%%%%%%%%%%%%%%%%%%%%%%%%%%%%%%%%%
\begin{Lem}
\label{locLip1}
Let $r \in (2,\infty)$, then there exists a constant $d_r>0$ depending only on $r$ such that  
the following inequality holds.
%%%%%%%%%%%%%%%%%%
\begin{equation}\label{locLip2}
   \begin{aligned}
    & \left|\left(|U|^{r - 2} u_i - |V|^{r - 2} v_i\right)(u_j - v_j)\right| 
       \leq d_r\left(|U|^{r - 2} + |V|^{r - 2}\right) |U - V|^2 
\\[2mm]
      & \qquad \qquad \qquad \text{ for all } \ U = (u_1, u_2), \  V = (v_1, v_2) \ \ \text{and} \ \  
           i, j = 1, 2. 
   \end{aligned}
\end{equation}
\end{Lem}
%%%%%%%%%%%%%%%%%%%%%%%%%%%%%%%%%%%%%%%%%%%%%%%%
%%%%%  Proof of Lemma 5     %%%%%%%%%%%%%%%%%%%%
%%%%%%%%%%%%%%%%%%%%%%%%%%%%%%%%%%%%%%%%%%%%%%%%
\begin{Prf} 
When $ |U| ~\! |V| = 0$, it is obvious that \eqref{locLip2} holds true with 
$d_r=1$. Then, in what follows, it suffices to consider the case where  $ |U| ~\! |V| \neq 0$.  
  We first note 
\begin{equation}\label{UgeqV}
  \begin{aligned}
   & \left| \left(|U|^{r - 2} u_i - |V|^{r - 2} v_i\right)(u_j - v_j) \right|
\\[2mm]
  = &  \   \left|  \left\{|U|^{r - 2} (u_i - v_i) 
        + \left(|U|^{r - 2} - |V|^{r - 2}\right) v_i\right\}(u_j - v_j) \right|
\\[2mm]
   \leq & \  |U|^{r - 2} |U - V|^2 
           +  \left| |U|^{r - 2} - |V|^{r - 2}  \right|
                   ~\! |V| ~\! |U - V|.
  \end{aligned}
\end{equation} 
Interchanging the roles of $U$ and $V$, we also get
\begin{equation}\label{VgeqU}
  \begin{aligned}
   & \left | \left(|U|^{r - 2} u_i - |V|^{r - 2} v_i\right)(u_j - v_j) \right|
\\[2mm]
   \leq & \ |V|^{r - 2} |V - U|^2 +  
            \left| |V|^{r - 2} - |U|^{r - 2}  \right| ~\! |U|  ~\! |V - U|.
  \end{aligned}
\end{equation}
  Here we claim that the following inequality holds. 
\begin{equation}\label{locLip3}
  \begin{aligned}
   \left| |U|^{r - 2} - |V|^{r - 2}\right| 
       \leq &  \ \tilde{d}_r\left(|U|^{r - 3} + |V|^{r - 3}\right) |U - V|,   
 \\[2mm] 
     &  \  \tilde{d}_r =
                 \begin{cases}
                   \frac{r -2}{2} & 2 < r \leq 3, \ \ 4 \leq r, 
\\[2mm]
                    \  1            & 3 < r< 4.
                 \end{cases}
  \end{aligned}     
\end{equation}
   In fact, we have
\begin{equation}\label{est:diff:UV}
  \begin{aligned}
     \left| |U|^{r - 2} - |V|^{r - 2} \right| 
       & = \left| \int_0^1 \frac{d}{d\theta} 
               \left\{ |V| + \theta ~\! (|U| - |V|)\right\}^{r - 2} ~\! d\theta \right| 
\\[2mm]
       & \leq (r - 2) \int_0^1 \left\{ |V| + \theta ~\! (|U| - |V|) \right\}^{r - 3} ~\! d\theta
                              \ |U - V|
\\[2mm]
       & =    (r - 2) \int_0^1 \left \{ \theta ~\! |U| 
                              + (1 - \theta) |V| \right\}^{r - 3} ~\! d\theta  \  |U - V|.
  \end{aligned}
\end{equation}
When \(2 < r \leq 3\) or \(4 \leq r\), by the convexity of the function 
\( s \mapsto |s|^{r - 3}\), we obtain 
\[
  \begin{aligned}
   \int_0^1 \left\{ \theta ~\! |U| + (1-  \theta) |V| \right\}^{r - 3} d\theta 
     & \leq \int_0^1 \theta ~\! |U|^{r - 3} + (1 - \theta) |V|^{r - 3} ~\! d\theta
\\[2mm]
       &  \leq \frac{1}{2}\left(|U|^{r - 3} + |V|^{r - 3}\right),
  \end{aligned}
\]
which together with \eqref{est:diff:UV} implies \eqref{locLip3}.
%%%%%%%%%%%%%%%%%%%%%%%%%%%%%%
As for the case \(3 < r < 4\), we note that the following inequality holds.
\begin{equation}\label{ineq:r34}
   (x + a)^{r - 3} \leq x^{r - 3} + a^{r - 3} \quad \text{for all} \  x>0, \ a>0.
\end{equation}
Indeed, put \(f(x) := (x + a)^{r - 3} - x^{r - 3} - a^{r - 3}\), then we get 
\begin{equation*} 
  f(0) = 0, \quad \text{and} \quad 
     f'(x) = (r - 3)\left\{(x + a)^{r - 4} - x^{r - 4}\right\} <0 
       \quad \text{for all} \ x>0, \ a>0.
\end{equation*}
Hence, applying \eqref{ineq:r34} with $x=\theta ~\! |U|$ and $a = (1-\theta) |V|$,
 we find 
\begin{equation*}
 \begin{aligned}
   \int_0^1 \left\{ \theta ~\! |U| + (1-  \theta) |V| \right\}^{r - 3} d\theta 
     & \leq \int_0^1 \theta^{r-3} ~\! |U|^{r - 3} + (1 - \theta)^{r-3} |V|^{r - 3} ~\! d\theta
\\[2mm]
       &  \leq \frac{1}{r-2}\left(|U|^{r - 3} + |V|^{r - 3}\right),
  \end{aligned}
\end{equation*}
which together with  \eqref{est:diff:UV} implies \eqref{locLip3}.

Combining (\ref{locLip3}) with (\ref{UgeqV}) and \eqref{VgeqU}, we get 
\begin{align}
     & \left|\left(|U|^{r - 2} u_i - |V|^{r - 2} v_i\right)(u_j - v_j)\right| 
          \leq \left\{|U|^{r - 2} + \tilde{d}_r\left(|U|^{r - 3}|V| + |V|^{r - 2}\right)\right\}|U - V|^2,
   \label{UgeqV1}
\\[2mm]
     & \left|\left(|U|^{r - 2} u_i - |V|^{r - 2} v_i\right)(u_j - v_j)\right| 
          \leq \left\{|V|^{r - 2} + \tilde{d}_r\left(|V|^{r - 3}|U| + |U|^{r - 2}\right)\right\}|U - V|^2.
   \label{VgeqU1}     
\end{align}
%%%%%%%%%%%%%%%%%%%%%%%%%%%%%%%%%%%
%%%%%  \( 4 , r\)       %%%%%%%%%%%
%%%%%%%%%%%%%%%%%%%%%%%%%%%%%%%%%%%
For the case where \(4 \leq r\), we note that Young's inequality gives 
\begin{equation}\label{est:UV:geq4}
        |U|^{r - 3}|V| \leq \frac{r - 3}{r - 2}|U|^{r - 2} + \frac{1}{r - 2}|V|^{r - 2},
\end{equation}
and plug this in \eqref{UgeqV1}, then we obtain (\ref{locLip2}) with 
$d_r = \displaystyle\frac{r-1}{2}$ . 
%%%%%%%%%%%%%%%%%%%%%%%%%%%%%%%%%%%
%%%%%  \(3 < r < 4\)    %%%%%%%%%%%
%%%%%%%%%%%%%%%%%%%%%%%%%%%%%%%%%%%

  As for the case where \(3 < r < 4\), interchanging the roles of $U$ and $V$ in 
   \eqref{est:UV:geq4} and adding the result to \eqref{est:UV:geq4}, we have 
\begin{equation}\label{est:UV:3geq4}    %###2017/09/09 add *
        |U|^{r - 3}|V| + |V|^{r-3} |U| 
            \leq |U|^{r - 2} + |V|^{r - 2}.
\end{equation}
  Then plugging \eqref{est:UV:3geq4} in the sum of \eqref{UgeqV1} and \eqref{VgeqU1}, 
   we derive (\ref{locLip2}) with $d_r = \displaystyle\frac{3}{2}$ .    %###2017/09/09 correct spelling
%   we derive (\ref{locLip2}) with $d_r = \dispalystyle\frac{3}{2}$ .
%%%%%%%%%%%%%%%%%%%%%%%%%%%%%%%%%%
%%%%%%   \(2 < r \leq 3\)     %%%%
%%%%%%%%%%%%%%%%%%%%%%%%%%%%%%%%%%

For the case where \(2 < r \leq 3\), we distinguish between two cases, namely  
$|U| \geq |V|$ or $|V| \geq |U|$. Suppose that $|U| \geq |V|$, 
then \(|U|^{r - 3} \leq |V|^{r - 3}\) holds, so we get by \eqref{UgeqV1}
\begin{equation}\label{UgeqV2}
   \left|\left(|U|^{r - 2} u_i - |V|^{r - 2} v_i\right)(u_j - v_j)\right| 
       \leq \left(|U|^{r - 2} + |V|^{r - 2}\right) |U - V|^2, 
\end{equation}
which gives \eqref{locLip2} with $ d_r = 1$.
For the case where \(|V| \geq |U|\), we repeat the same argument as above 
with \eqref{UgeqV1} replaced by \eqref{VgeqU1} to obtain (\ref{locLip2}) with $d_r=1$.
\end{Prf}
%%%%%%%%%%%%%%%%%%%%%%%%%%%%%%%%%%%%%%%%%%%%%%%%%%%%%%%%%%
%%%%%%     Corollary 5      %%%%%%%%%%%%%%%%%%%%%%%%%%%%%%
%%%%%%%%%%%%%%%%%%%%%%%%%%%%%%%%%%%%%%%%%%%%%%%%%%%%%%%%%%
   The following estimate follows directly from Lemma \ref{locLip1}.
\begin{Cor}
\label{locLip4}
There exists a constant \(C\) such that the following estimates hold 
for all $U, \ V \in D(\partial \psi_r)$. 
\begin{align}
\label{ppsiUppsiVUV}
   \left|(\partial\psi_r(U) - \partial\psi_r(V), U - V)_{\mathbb{L}^2}\right| 
      &  \leq C\left(\psi_r(U)^{r - 2} + \psi_r(V)^{r - 2}\right)|U - V|_{\mathbb{L}^r}^2,
\\[2mm]
\label{ppsiUppsiVIUV}
    \left|(\partial\psi_r(U) - \partial\psi_r(V), I(U - V))_{\mathbb{L}^2}\right| 
         & \leq C\left(\psi_r(U)^{r - 2} + \psi_r(V)^{r - 2}\right)|U - V|_{\mathbb{L}^r}^2.
\end{align}
\end{Cor}
%%%%%%%%%%%%%%%%%%%%%%%%%%%%%%%%%%%%%%%%%%%%%%%%%%%%%%%%
%%%%%    Proof of Uniqueness      %%%%%%%%%%%%%%%%%%%%%%
%%%%%%%%%%%%%%%%%%%%%%%%%%%%%%%%%%%%%%%%%%%%%%%%%%%%%%%%
We now proceed to the proof of the uniqueness.
\begin{Prf}[Theorem \ref{lwpndbd} (Uniqueness)]
Let \(U, V\) be two solutions of (ACGL) with \(U(0) = U_0\) and \(V(0) = V_0\) on \([0, T_0]\).
Multiplying the difference of two equations by \( W := U - V\) and using the linearity of \(\partial\varphi\), (\ref{ortho}) and Corollary \ref{locLip4} with $r=q$, we get
\begin{equation}\label{DUV}
\begin{aligned}
   \frac{1}{2}\frac{d}{dt}|W|_{\mathbb{L}^2}^2 + 2 \lambda ~\! \varphi (W) 
      & \leq \gamma_+ |W|_{\mathbb{L}^2}^2 
          + \left((\kappa + I\beta)(\partial\psi_q(U) - \partial\psi_q(V)), W \right)_{\mathbb{L}^2}
\\[2mm]
      & \leq \gamma_+ |W|_{\mathbb{L}^2}^2 
              + C \left( \psi_q(U)^{q - 2} + \psi_q(V)^{q - 2} \right) |W|_{\mathbb{L}^q}^2,
\end{aligned}
\end{equation}
where \(C\) is a constant depending only on \(q, \kappa, \beta\).

By our assumption on \(q\) being Sobolev subcritical, using the parameter \(\eta \in (0,1)\) defined by 
\begin{equation}\label{eta}
  \frac{1}{q} = \left(\frac{1}{2} - \frac{1}{N}\right)(1 - \eta) + \frac{\eta}{2},
\end{equation}
we obtain
\begin{equation}
\label{GNH1-1}
|W|_{\mathbb{L}^q} \leq C \left( \varphi(W) \right)^{\frac{1 - \eta}{2}}|W|_{\mathbb{L}^2}^{\eta},
\end{equation}
with an appropriate constant \(C\).
   Thus by (\ref{DUV}), (\ref{GNH1-1}) and Young's inequality, we obtain 
\begin{equation}\label{DiffUV}
   \frac{1}{2}\frac{d}{dt}|W|_{\mathbb{L}^2}^2 
     + \lambda ~\! \varphi(W) 
        \leq \left[ C \left(\psi_q(U)^{q - 2} + \psi_q(V)^{q - 2}\right)^{\frac{1}{\eta}}
                   + \gamma_+ \right] |W|_{\mathbb{L}^2}^2,
\end{equation}
where \(C\) depends on \(\lambda, \kappa, \beta, \gamma, \eta\).
%%%%%
Since the regularity (i) and (iii) in Theorem 1 assures the absolute continuity of  \(\varphi(U)\) and \(\varphi(V)\) on \([0, T_0]\)
(see \cite{B1}), we see that \(\psi_q(U)\) and \(\psi_q(V)\) are uniformly bounded above by a positive constant \(M\) on \([0, T_0]\).
Then applying Gronwall's inequality to (\ref{DiffUV}), we obtain
\begin{equation}\label{UNI}
   |U(t) - V(t)|_{\mathbb{L}^2}^2 
      \leq |U_0 - V_0|_{\mathbb{L}^2}^2 ~\! 
              e^{[2 C (2 M^{q - 2})^{\frac{1}{\eta}} + \gamma_+] ~\! t} \quad 
                  \forall t \in [0,T_0],
\end{equation}
whence follows the uniqueness.
\end{Prf}

%%%%%%%%%%%%%%%%%%%%%%%%%%%%%%%%%%%%%%%%%%%%%%%%%%
%%%%%%%%%%%%%%%%%%%%%%%%%%%%%%%%%%%%%%%%%%%%%%%%%%
%%%%%%%%%%   Pf of Thm 3   %%%%%%%%%%%%%%%%%%%%%%%
%%%%%%%%%%%%%%%%%%%%%%%%%%%%%%%%%%%%%%%%%%%%%%%%%%
%%%%%%%%%%%%%%%%%%%%%%%%%%%%%%%%%%%%%%%%%%%%%%%%%%
\section{Proof of Theorem \ref{gebd}}
First we prepare some lemmas.
%%%%%%%%%%%%%%%%%%%%%%%%%%%%%%%%%%%%%%%%%%%%%%%%%
%%%%%%    Lemma 6    %%%%%%%%%%%%%%%%%%%%%%%%%%%%
%%%%%%%%%%%%%%%%%%%%%%%%%%%%%%%%%%%%%%%%%%%%%%%%%
\begin{Lem}\label{coer}
Let all the assumptions in Theorem \ref{gebd} be satisfied.
There exist \(\varepsilon_0 > 0\) and \(\delta > 0\) such that for all \(U \in D(\varphi) = \mathbb{H}_0^1(\Omega)\) satisfying \(\varphi(U) < \varepsilon_0\), it holds that
\begin{equation}\label{coer1}
    (\lambda ~\! \partial \varphi( U ) - \kappa ~\! \partial \psi_q (U) 
           - \gamma ~\!  U, U)_{\mathbb{L}^2}    
                 \geq \delta ~\! \varphi(U) 
                    \geq \frac{\delta\lambda_1}{2} |U|_{\mathbb{L}^2}^2.
\end{equation}
\end{Lem}
\begin{Prf}
%%%%%%%%%%%%%%%%%%%%%%%%%%%%%%%%%%%%%%%%%%%%%%%%
%%%%%    Proof of Lemma 6      %%%%%%%%%%%%%%%%%
%%%%%%%%%%%%%%%%%%%%%%%%%%%%%%%%%%%%%%%%%%%%%%%%
Since we assume \(\lambda ~\! \lambda_1 > \gamma\), we first note that
\begin{equation}\label{7.2}
   (\lambda ~\! \partial\varphi(U) - \gamma ~\! U, U)_{\mathbb{L}^2} 
     =  2 ~\! \lambda ~\! \varphi(U) - \gamma ~\! |U|_{\mathbb{L}^2}^2 
        \geq \delta_0 ~\! \varphi(U), \ \delta_0 := 2\left(\lambda - \frac{\gamma}{\lambda_1}\right) > 0.
\end{equation}
On the other hand, since \(q\) is subcritical, there exist a constant \(C > 1\) such that
\[
\psi_q(U) \leq C\varphi(U)^{\frac{q}{2}}\quad \forall U \in \mathbb{H}^1_0(\Omega).
\]
Hence we have
\begin{equation}\label{7.3}
      ( - \kappa ~\! \partial\psi_q(U), U)_{\mathbb{L}^2} 
         = - \kappa ~\! q ~\! \psi_q(U) 
           \geq - C ~\!\! \kappa ~\!  q ~\!  (\varphi(U))^{\frac{q}{2}} 
                \quad \forall U \in \mathbb{H}^1_0(\Omega).
\end{equation}
Here we take
\begin{equation}\label{7.4}
     \delta = \frac{\delta_0}{2} > 0, \quad 
        \varepsilon_0 = \left( \frac{\delta_0}{2 ~\! C ~\! \kappa q} \right)^{\frac{2}{q - 2}} > 0.
\end{equation}
Then, by (\ref{7.2}), (\ref{7.3}) and (\ref{7.4}), we see that
\[
   \begin{aligned}
     (\lambda ~\! \partial \varphi(U) - \kappa ~\! \partial \psi_q(U) 
                                - \gamma ~\! U, U)_{\mathbb{L}^2} 
        & \geq \delta_0 \varphi(U) - C ~\!\! \kappa ~\!\! q(\varphi(U))^{\frac{q}{2}}
\\[2mm]
        & \geq \frac{\delta_0}{2}\varphi(U)
   \end{aligned}
\]
holds if \(\varphi(U) < \varepsilon_0\).
\end{Prf}
%%%%%%%%%%%%%%%%%%%%%%%%%%%%%%%%%%%%%%%%%%%%%%%%%%%%%%
%%%%%%   Lemma 7    %%%%%%%%%%%%%%%%%%%%%%%%%%%%%%%%%%
%%%%%%%%%%%%%%%%%%%%%%%%%%%%%%%%%%%%%%%%%%%%%%%%%%%%%%
\begin{Lem}
\label{Gtyineq}
   Let \(f(t) \in {\rm L}^1(0, T)\) and \(j(t)\) be an absolutely 
     continuous positive function on \([0, S]\) with \(0 < S \leq T\) such that
\begin{equation}\label{7.5}
    \frac{d}{dt} ~\! j(t) + \delta~\! j(t) \leq K ~\! |f(t)|   \quad \mbox{a.e.}\ t \in [0, S],
\end{equation}
where \(\delta > 0\) and \(K > 0\). Then we have
\begin{equation}\label{7.6}
  \begin{aligned}
    & j(t) \leq j(0) ~\! e^{-\delta t} 
       + \frac{K}{1 - e^{-\delta}} ~\! \opnorm{f}_1 
           \quad  \forall t \in [0, S],
\\[2mm]
    & \opnorm{f}_1 = \sup\left\{\int_S^{S + 1}|\tilde{f}(t)| ~\! dt 
        \,;\, 0 \leq S < \infty\right\},
\end{aligned}
\end{equation}
where \(\tilde{f}\) is the zero extension of \(f\) to \([0, \infty)\).
\end{Lem}
%%%%%%%%%%%%%%%%%%%%%%%%%%%%%%%%%%%%%%%%%%%%%%%%%%%
%%%%%%%%%%%   Proof of Lemma 7      %%%%%%%%%%%%%%%
%%%%%%%%%%%%%%%%%%%%%%%%%%%%%%%%%%%%%%%%%%%%%%%%%%%
\begin{Prf}
This fact is essentially proved in Lemma 4.3 of \^Otani \cite{O2}. 
  Let $ n = [t]$, i.e., $n \in \mathbb{N}\cup \{0\}$ such that 
   $ n \leq t < n+1$. Then it is easy to see that (\ref{7.5}) implies
\[
  \begin{aligned}
     j(t) 
       & \leq j(0) ~\! e^{-\delta t} + \int_0^t K ~\! |f(s)|e^{-\delta(t - s)} ~\! ds
\\
       & \leq j(0) ~\! e^{-\delta t} + K \sum_{k=0}^{n-1} \int_{t-k-1}^{t-k} | \tilde{f}(s)|e^{-\delta(t - s)} ~\! ds 
                     + K \int_0^{t-n}| \tilde{f}(s)|e^{-\delta(t - s)} ~\! ds 
\\
       & \leq j(0) ~\! e^{-\delta t} +  K \sum_{k=0}^{n-1} \int_{t-k-1}^{t-k} | \tilde{f}(s)|e^{- k \delta} ~\! ds 
                     + K \int_0^{t-n}| \tilde{f}(s)|e^{- n \delta} ~\! ds 
\\
       & \leq j(0) ~\! e^{-\delta t} + \frac{K}{1 - e^{-\delta}} ~\! \opnorm{f}_1.
  \end{aligned}
\]
\end{Prf}
%%%%%%%%%%%%%%%%%%%%%%%%%%%%%%%%%%%%%%%%%%%%%%%%%%%%%%%%%%%%
%%%%%%%%%%%   Lemma 8      %%%%%%%%%%%%%%%%%%%%%%%%%%%%%%%%%
%%%%%%%%%%%%%%%%%%%%%%%%%%%%%%%%%%%%%%%%%%%%%%%%%%%%%%%%%%%%
\begin{Lem}
\label{gbdd}
Let all the assumptions in Theorem \ref{gebd} be satisfied.
Then there exist \(\varepsilon_1 > 0\) and \(N > 0\) independent of \(T\) such that for any \(r \in (0, \varepsilon_1)\), if \(\varphi(U_0) \leq r^2\) and \(\opnorm{F}_2 \leq r\), then any solution \(U\) of (ACGL) on \([0, S] \ (S \in (0,T]) \) satisfying (i)-(iii) of Theorem \ref{lwpndbd} satisfies
\begin{equation}
\label{gbdd1}
\varphi(U(t)) < Nr^2\quad \forall t \in [0, S].
\end{equation}
\end{Lem}
%%%%%%%%%%%%%%%%%%%%%%%%%%%%%%%%%%%%%%%%%%%%%%%%%%%%%%%%%%%%
%%%%%%%%%   Proof of Lemma 8       %%%%%%%%%%%%%%%%%%%%%%%%%
%%%%%%%%%%%%%%%%%%%%%%%%%%%%%%%%%%%%%%%%%%%%%%%%%%%%%%%%%%%%
\begin{Prf}
By the same argument as that for (\ref{AEhpphiUh}), we get
\begin{equation}\label{7.8}
       \frac{d}{dt}\varphi(U) + \frac{\lambda}{4} ~\! |\partial\varphi(U)|_{\mathbb{L}^2}^2 
          \leq 2 ~\! \gamma_+ \varphi(U) 
                   + \frac{\kappa^2 + \beta^2}{\lambda}
                        ~\! |\partial\psi_q(U)|_{\mathbb{L}^2}^2 
                           + \frac{1}{\lambda}~\! |F|_{\mathbb{L}^2}^2.
\end{equation}
Then using (\ref{ppsiqUlL}) with \(\varepsilon = \varepsilon_\lambda := \lambda^2/8(\kappa^2 + \beta^2)\), we obtain
\begin{equation}\label{7.9}
    \frac{d}{dt} \varphi(U) + \frac{\lambda}{8} ~\! |\partial\varphi(U)|_{\mathbb{L}^2}^2 
       \leq 2 ~\! \gamma_+ \varphi(U) + C_{\varepsilon_\lambda}\varphi(U)^\chi 
           + \frac{1}{\lambda} ~\! |F|_{\mathbb{L}^2}^2,
\end{equation}
where \(\chi= \xi(q - 1)/\{1 - (q - 1)(1 - \xi)\} > 1\).

Here we fix \(N\) and \(\varepsilon_1\) by
\begin{align}
    N  & = \left[2 + \frac{1}{1 - e^{-\frac{\lambda\lambda_1}{4}}}\left\{(2\gamma_+C_\lambda)N_2 
            + \frac{1}{\lambda}\right\}\right],
\\
\notag 
    N_2 & = \frac{1}{\delta} \left( N_1 + \frac{1}{2} N_1^2 \right), \quad N_1 
           = \sqrt{\frac{2}{\lambda_1}} + \frac{1}{1 - e^{- \frac{\delta\lambda_1}{2}}}.
\\[1mm]
\label{7.11}
   \varepsilon_1 & = \frac{\varepsilon_0}{N} 
          \quad  \mbox{(\(\varepsilon_0\) is the number appearing in Lemma \ref{coer}).}
\end{align}
Then we claim that (\ref{gbdd1}) holds true for all \(t \in [0, S]\).
Suppose that this is not the case, then by the continuity of \(\varphi(U(t))\), there exists \(t_1 \in (0, S)\) such that
\begin{equation}
       \varphi(U(t)) < Nr \quad   \forall t \in [0, t_1)  
            \ \  \mbox{and} \ \ \varphi(U(t_1)) = N r.
\end{equation}
We are going to show that this leads to a contradiction.
 We first multiply (ACGL) by \(U(t)\) for \(t \in [0, t_1]\).
Then since \(\varphi(U(t)) \leq Nr \leq \varepsilon_0\) for all \(t \in [0, t_1]\), 
  Lemma \ref{coer} gives
\begin{align}\label{7.13}
    & \frac{1}{2}\frac{d}{dt}|U(t)|_{\mathbb{L}^2}^2 
        + \delta\varphi(U(t)) \leq |F(t)|_{\mathbb{L}^2}|U(t)|_{\mathbb{L}^2} 
          & \forall t \in [0, t_1], 
\\[2mm]
\label{7.14}
    &  \frac{d}{dt}|U(t)|_{\mathbb{L}^2} 
        + \frac{\delta\lambda_1}{2}|U(t)|_{\mathbb{L}^2} \leq |F(t)|_{\mathbb{L}^2} 
          & \forall t \in [0, t_1].
\end{align}
Hence, by (\ref{7.14}) and Lemma \ref{Gtyineq}, we get
\begin{equation}
    \sup_{0 \leq t \leq t_1}|U(t)|_{\mathbb{L}^2} 
         \leq \left(\sqrt{\frac{2}{\lambda_1}} 
                 + \frac{1}{1 - e^{-\frac{\delta\lambda_1}{2}}}\right)r = N_1 r,
\end{equation}
where we used the fact that 
  \(|U(0)|_{\mathbb{L}^2} \leq \sqrt{\frac{2}{\lambda}\varphi(U(0))} 
                             \leq \sqrt{\frac{2}{\lambda}}r\) 
                                and \( \opnorm{|F(t)|_{\mathbb{L}^2}}_1 = \opnorm{F}_2 \leq r\).
Hence the integration of (\ref{7.13}) over \((t, t+1)\) gives
\begin{equation}\label{7.16}
    \sup_{0 \leq t < \infty} \int_t^{t + 1} \tilde{\varphi}(U(\tau))d\tau 
              \leq \frac{1}{\delta} \left( N_1 r^2 + \frac{1}{2} ~\! N_1^2 ~\! r^2 \right) = N_2 ~\! r^2,
\end{equation}
where \(\tilde{\varphi}(U(\cdot))\) is the zero extension of \(\varphi(U(\cdot))\) to \([0, \infty)\).
Here we note
\[
|\partial\varphi(U)|_{\mathbb{L}^2}|U|_{\mathbb{L}^2} \geq |(\partial\varphi(U), U)_{\mathbb{L}^2}| = 2\varphi(U) \geq 2\sqrt{\varphi(U)}\sqrt{\frac{\lambda_1}{2}}|U|_{\mathbb{L}^2},
\]
whence follows
\begin{equation}
|\partial\varphi(U)|_{\mathbb{L}^2}^2 \geq 2\lambda_1\varphi(U)\quad\forall U \in {\rm D}(\partial\varphi),
\end{equation}
which together with (\ref{7.9}) yields
\begin{equation}\label{7.18}
    \frac{d}{dt} \varphi(U) + \frac{\lambda ~\! \lambda_1}{4} ~\! \varphi(U) 
        \leq 2 ~\! \gamma_+\varphi(U) + C_{\varepsilon_\lambda} \varphi(U)^\chi 
         + \frac{1}{\lambda}~\! |F|_{\mathbb{L}^2}^2 
              \quad \forall t \in [0, t_1].
\end{equation}
Without loss of generality, we can take \(Nr \leq N\varepsilon_1 \leq \varepsilon_0 \leq 1\).
Then Since \(\chi > 1\), in view of (\ref{7.16}), (\ref{7.18}) and Lemma \ref{Gtyineq}, we obtain
\begin{equation}
\varphi(U(t)) \leq \left[1 + \frac{1}{1 - e^{-\frac{\lambda\lambda_1}{4}}}\left\{(2\gamma_+ + C_{\varepsilon_\lambda}N_2 + \frac{1}{\lambda}\right\}\right]r^2 < Nr^2,
\end{equation}
which contradicts (\ref{7.11}).
\end{Prf}

%%%%%%%%%%%%%%%%%%%%%%%%%%%%%%%%%%%%%%%%%%%%%%%%%%
%%%%%%%%%%   Pf of Thm 3   %%%%%%%%%%%%%%%%%%%%%%%
%%%%%%%%%%%%%%%%%%%%%%%%%%%%%%%%%%%%%%%%%%%%%%%%%%
Now we are ready to prove Theorem \ref{gebd}.
\begin{Prf}[Theorem \ref{gebd}]
By virtue of Theorem \ref{altbd}, in order to show the existence of global solution, it suffices to give the a priori bound for \(\varphi(U(t))\).
Indeed due to Lemma \ref{gbdd}, we obtain the a priori bound (\ref{gbdd1}) for \(\varphi(U(t))\).
\end{Prf}

%%%%%%%%%%%%%%%%%%%%%%%%%%%%%%%%%%%%%%%%%%%%%%%%%%
%%%%%%%%%%%%%%%%%%%%%%%%%%%%%%%%%%%%%%%%%%%%%%%%%%
%%%%%%%%%%   Bibliographies   %%%%%%%%%%%%%%%%%%%%
%%%%%%%%%%%%%%%%%%%%%%%%%%%%%%%%%%%%%%%%%%%%%%%%%%
%%%%%%%%%%%%%%%%%%%%%%%%%%%%%%%%%%%%%%%%%%%%%%%%%%

%%%%%%%%%%%%%%%%%%%%%%%%%%%%%%%%%%%%%%%%%%%%%%%%%%
%%%%%%%%%%%%%%%%%%%%%%%%%%%%%%%%%%%%%%%%%%%%%%%%%%
%%%%%%%%%%%%%%%%%%%%%%%%%%%%%%%%%%%%%%%%%%%%%%%%%%
%%%%%%%%%%%%%%%%%%%%%%%%%%%%%%%%%%%%%%%%%%%%%%%%%%
%%%%%%%%%%%%%%%%%%%%%%%%%%%%%%%%%%%%%%%%%%%%%%%%%%

\section*{References}
\begin{thebibliography}{99}
%
\bibitem%[(1973)]
{B1} Br\'ezis, H., 
{\it Op\'erateurs Maximaux Monotones et Semi-Groupes de Contractions dans les Espaces de Hilbert}, 
North-Holland (1973), \\
(https://doi.org/10.1016/s0304-0208(08)x7125-7.)
%
\bibitem%[(1971)]
{B2} Br\'ezis, H., 
``Monotonicity methods in Hilbert spaces and some applications to nonlinear partial differential equations,'' 
in E. H. Zarantonello (Ed.), {\it Contributions to Nonlinear Functional Analysis (Proc. Sympos., Math. Res. Center, Univ. Wisconsin, Madison, Wis., 1971)}, Academic Press, New York (1971) 101-156, \\
(https://doi.org/10.1016/B978-0-12-775850-3.50009-1).
%
\bibitem%[(2014)]
{CDF1} Cazenave, T.; Dias, J.; Figueira, M., 
``Finite-time blowup for a complex Ginzburg-Landau equation with linear driving,'' 
{\it J. Evol. Equ.} {\bf 14} (2014), no. 2, 403-415, \\
(https://dx.doi.org/10.1007/s00028-014-0220-z).
%
\bibitem%[(2013)]
{CDW1} Cazenave, T.; Dickstein, F.; Weissler, B., 
``Finite-Time Blowup for a complex Ginzburg-Landau equation,'' 
{\it SIAM J. Math. Anal.} {\bf 45} (2013), no. 1, 244-266, \\
(https://dx.doi.org/10.1137/120878690).
%
\bibitem%[(2014)]
{CDW2} Cazenave, T.; Dickstein, F.; Weissler, F. B., 
``Standing waves of the complex Ginzburg-Landau equation,'' 
{\it Nonlinear Anal.} {\bf 103} (2014), 26-32, \\
(https://dx.doi.org/10.1016/j.na.2014.03.001).
%
\bibitem%[(1993)]
{CH1} Cross, C.; Hohenberg, C., 
``Pattern formation outside of equilibrium,'' 
{\it Rev. Mod. Phys.} {\bf 65} (1993), 851-1112, \\
(https://dx.doi.org/10.1103/RevModPhys.65.851).
%
\bibitem%[(1950)]
{GL1} Ginzburg, L.; Landau, D., 
``On the theory of superconductivity,''
{\it Zh. Eksp. Teor. Fiz.} {\bf 20} (1950), 1064-1082 (in Russian); 
translation in D. ter Haar (Ed.), {\it Collected Papers of L.D. Landau}, Pergamon Press (1965), 546-568, \\
(https://doi.org/10.1016/b978-0-08-010586-4.50078-x).
%
\bibitem%[(accepted)]
{KOS1} Kuroda, T.; \^Otani, M.; Shimizu, S., 
``Initial-boundary value problems for complex Ginzburg-Landau equations in general domains,'' 
{\it Adv. Appl. Math. Sci.} {\bf 26} (2017), 119-141.
%
\bibitem%[(2009)]
{N1} Nishiura, Y., ``Far-from-Equilibrium Dynamics,''
{\it Translations of Mathematical Monographs} {\bf 209}, AMS (2002).%% Iwanami Series in Modern Mathematics
%
\bibitem%[(1977)]
{O1} \^Otani, M., 
``On the existence of strong solutions for \(\frac{du}{dt}(t) + \partial\psi_1(u(t)) - \partial\psi_2(u(t)) \ni f(t)\),'' 
{\it J. Fac. Sci. Univ. Tokyo Sect. IA Math.} {\bf 24} (1977), no. 3, 575-605, \\
(https://hdl.handle.net/2261/6206).
%
\bibitem%[(1982)]
{O2} \^Otani, M., 
``Nonmonotone perturbations for nonlinear parabolic equations associated with subdifferential operators, Cauchy Problems,'' 
{\it J. Differential Equations} {\bf 46} (1982), no. 2, 268-299, \\
(https://dx.doi.org/10.1016/0022-0396(82)90119-X).
%
\bibitem%[(2002)]
{OY1} Okazawa, N.; Yokota, T., 
``Global existence and smoothing effect for the complex Ginzburg-Landau equation with p-Laplacian,'' 
{\it J. Differential Equations} {\bf 182} (2002), no. 2, 541-576, \\
(https://dx.doi.org/10.1006/jdeq.2001.4097).
%
\bibitem%[(2016)]
{SYY2} Shimotsuma, D.; Yokota, T.; Yoshii, K., 
``Existence and decay estimates of solutions to complex Ginzburg-Landau type equations,'' 
{\it J. Differential Equations} {\bf 260} (2016), no. 3, 3119-3149, \\
(https://dx.doi.org/10.1016/j.jde.2015.10.030).
%
\bibitem%[(1988)]
{T1} Temam, R., 
{\it Infinite-dimensional dynamical systems in mechanics and physics}, 
Applied Mathematical Sciences, {\bf 68}, Springer-Verlag, New York, (1988), \\
(https://doi.org/10.1007/978-1-4684-0313-8).
%
\bibitem%[(1990)]
{Y1} Yang, Y.,
``On the Ginzburg-Landau wave equation,''
{\it Bull. London Math. Soc.} {\bf 22} (1990), no. 2, 167-170, \\
(https://dx.doi.org/10.1112/blms/22.2.167).
%
\end{thebibliography}
\end{document}